\documentclass[12pt,oneside,reqno]{amsart}
\usepackage{mathrsfs}
\usepackage{graphics}
\usepackage{cite}
\usepackage{amssymb}
\usepackage{enumerate}
\newcommand{\dif}{\mathrm{d}}

\newcommand{\be}{\begin{eqnarray}}
\newcommand{\ee}{\end{eqnarray}}
\newcommand{\ce}{\begin{eqnarray*}}
\newcommand{\de}{\end{eqnarray*}}
\newtheorem{theorem}{Theorem}[section]
\newtheorem{lemma}[theorem]{Lemma}
\newtheorem{remark}[theorem]{Remark}
\newtheorem{definition}[theorem]{Definition}
\newtheorem{proposition}[theorem]{Proposition}
\newtheorem{Example}[theorem]{Example}
\newtheorem{corollary}[theorem]{Corollary}
\newtheorem{condition}[theorem]{Condition}
\def\e{\varepsilon}
\def\s{\sigma}
\def\t{\theta}
\def\a{\alpha}

\def\b{\beta}
\def\d{\delta}
\def\p{\partial}
\def\g{\gamma}
\def\l{\lambda}

\def\[{{\Big[}}
\def\]{{\Big]}}
\def\<{{\langle}}
\def\>{{\rangle}}
\def\({{\Big(}}
\def\){{\Big)}}

\def\no{\nonumber}
\def\bt{\begin{theorem}}
\def\et{\end{theorem}}
\def\bl{\begin{lemma}}
\def\el{\end{lemma}}
\def\br{\begin{remark}}
\def\er{\end{remark}}
\def\bx{\begin{Example}}
\def\ex{\end{Example}}
\def\bd{\begin{definition}}
\def\ed{\end{definition}}
\def\bp{\begin{proposition}}
\def\ep{\end{proposition}}
\def\bc{\begin{corollary}}
\def\ec{\end{corollary}}
\def\bco{\begin{condition}}
\def\eco{\end{condition}}

\def\cD{{\mathcal D}}

\def\cQ{{\mathcal Q}}

\def\mE{{\mathbb E}}

\def\mM{{\mathbb M}}
\def\mN{{\mathbb N}}

\def\mP{{\mathbb P}}

\def\mR{{\mathbb R}}
\def\mS{{\mathbb S}}

\def\sA{{\mathscr A}}
\def\sB{{\mathscr B}}

\def\sL{{\mathscr L}}

\def\sV{{\mathscr V}}

\def\geq{\geqslant}
\def\leq{\leqslant}

\begin{document}

\allowdisplaybreaks

\title{Large deviation principles for fully coupled multiscale multivalued stochastic systems}

\author{Huijie Qiao}

\dedicatory{School of Mathematics,
Southeast University,\\
Nanjing, Jiangsu 211189, P.R.China\\
hjqiaogean@seu.edu.cn}



\thanks{This work was supported by NSF of China (No.12071071) and the Jiangsu Provincial Scientific Research Center of Applied Mathematics (No. BK20233002).}


\date{}

\begin{abstract}
This study focuses on large deviation principles for fully coupled multiscale multivalued stochastic systems, in which the slow component is governed by a multivalued stochastic differential equation and the fast component is described by a general stochastic differential equation. First, we establish the large deviation principle for the slow component at any fixed time by leveraging viscosity solutions of second-order Hamilton-Jacobi-Bellman equations involving multivalued operators. Subsequently, we illustrate the theoretical results through a concrete example.
\end{abstract}

\maketitle \rm

\section{Introduction}

In this paper, we investigate the large deviation principle (LDP for short) for the slow component of the following slow-fast system defined on $\mathbb{R}^{n} \times \mathbb{R}^{m}$: for any $T>0$,
\begin{equation}
\left\{
\begin{array}{l}
dX_{t}^{\varepsilon,\gamma} \in -A(X_{t}^{\varepsilon,\gamma})dt + b_{1}(X_{t}^{\varepsilon,\gamma},Y_{t}^{\varepsilon,\gamma})dt + \sqrt{\varepsilon}\sigma_{1}(X_{t}^{\varepsilon,\gamma},Y_{t}^{\varepsilon,\gamma})dW^1_{t}, \\
X_{0}^{\varepsilon,\gamma} = x_0 \in \overline{\mathcal{D}(A)}, \quad 0 \leq t \leq T, \\
dY_{t}^{\varepsilon,\gamma} = \frac{1}{\gamma}b_{2}(X_{t}^{\varepsilon,\gamma},Y_{t}^{\varepsilon,\gamma})dt + \frac{1}{\sqrt{\gamma}}\sigma_{2}(X_{t}^{\varepsilon,\gamma},Y_{t}^{\varepsilon,\gamma})dW^2_{t}, \\
Y_{0}^{\varepsilon,\gamma} = y_0, \quad 0 \leq t \leq T,
\end{array}
\right.
\label{Eq1i}
\end{equation}
where $(W^1_{t})$ and $(W^2_{t})$ are independent $d_1$- and $d_2$-dimensional standard Brownian motions defined on a complete filtered probability space $(\Omega,\mathcal{F},\{\mathcal{F}_{t}\}_{t\in[0,T]},\mathbb{P})$, respectively. The operator $A$ is a multivalued maximal monotone operator (see Subsection \ref{mmo}), and the mappings $b_{1}:\mathbb{R}^{n}\times\mathbb{R}^{m}\rightarrow\mathbb{R}^{n}$, $\sigma_{1}:\mathbb{R}^{n}\times\mathbb{R}^{m}\rightarrow\mathbb{R}^{n\times d_1}$, $b_{2}:\mathbb{R}^{n}\times\mathbb{R}^{m}\rightarrow\mathbb{R}^{m}$, and $\sigma_{2}:\mathbb{R}^{n}\times\mathbb{R}^{m}\rightarrow\mathbb{R}^{m\times d_2}$ are all Borel measurable. Here, $0 < \varepsilon < 1$ is a small parameter, and $0 < \gamma < 1$ represents another small parameter that characterizes the timescale separation between the processes $X_{\cdot}^{\varepsilon,\gamma}$ and $Y_{\cdot}^{\varepsilon,\gamma}$. Moreover, $\gamma$ is assumed to be a function of $\varepsilon$.

Systems of the form (\ref{Eq1i}) are commonly referred to as fully coupled multiscale stochastic systems and have found applications in various fields, including biology, chemistry, physics, and finance (\cite{fff, ffk, fk}). When $A=0$, numerous results concerning the LDP for system (\ref{Eq1i}) have been established (\cite{bz, fff, ffk, hlls, yK, ls, rl, ll, aP, kSp, aV1, aV2}). In the case where $A = \partial \mathcal{I}_\mathcal{O}$, with $\mathcal{O}$ being a closed convex subset of $\mathbb{R}^n$ and
$$
\mathcal{I}_{\mathcal{O}}(x):= 
\begin{cases}
0, & \text{if } x \in \mathcal{O}, \\ 
+\infty, & \text{if } x \notin \mathcal{O},
\end{cases}
$$
and under the assumptions that $\sigma_1$ is independent of $y$ and $\lim_{\varepsilon\to 0} \gamma/\varepsilon = 0$, Kushner \cite{ku} briefly discussed the validity of the LDP for the slow component in the path space. More recently, in \cite{q1}, we incorporated a multivalued maximal monotone operator into the fast component of system (\ref{Eq1i}), maintained the same assumptions on $\sigma_1$ and $\gamma/\varepsilon$, and rigorously established the LDP for the slow component in the path space using the weak convergence approach. Subsequently, in \cite{q2}, we extended this result to multiscale multivalued McKean-Vlasov stochastic systems. A natural question arises: does the LDP still hold for system (\ref{Eq1i}) when $\sigma_1$ depends on $y$ and $\lim_{\varepsilon\to 0} \gamma/\varepsilon = 0$?

Currently, three primary methods are employed to study the LDP for fully coupled multiscale stochastic systems: exponential tightness (\cite{aV1, aV2}), tightness of occupation measures (\cite{bz, ds, rl, aP, kSp}), and viscosity solutions of second-order Hamilton-Jacobi-Bellman equations (\cite{ffk, ls, kp}). However, due to the presence of the multivalued maximal monotone operator $A$ in system (\ref{Eq1i}), the first two approaches are not directly applicable. Consequently, we establish the LDP for system (\ref{Eq1i}) by means of viscosity solutions associated with a second-order Hamilton-Jacobi-Bellman equation involving the multivalued operator $A$.

The application of viscosity solutions of second-order Hamilton-Jacobi-Bellman equations to LDPs for stochastic differential equations (SDEs for short) was first introduced by Evans and Ishii \cite{ei}. Subsequently, Feng and Kurtz \cite{fk} developed a general framework based on this approach. To date, several LDP results have been obtained using this methodology (\cite{mK, rw, aS}).

The application of this method to LDPs for fully coupled multiscale stochastic systems was first proposed by Feng, Fouque, and Kumar in \cite{ffk}. Subsequently, Kumar and Popovic \cite{kp} employed this approach to establish the LDP for multiscale jump-diffusion processes, while Li and Shao \cite{ls} proved two LDPs for two-time-scale regime-switching processes. However, applying this method to derive the LDP for system (\ref{Eq1i}) is non-trivial. First, due to the presence of the multivalued maximal monotone operator $A$, a new viscosity solution framework must be developed for the associated second-order Hamilton-Jacobi-Bellman equation involving $A$. Second, the existence of $A$ prevents the use of supermartingale arguments when verifying the exponential tightness of $\{X^{\varepsilon,\gamma}_t, \varepsilon, \gamma > 0\}$ for any $t \in [0,T]$ (see Section 4.5 in \cite{fk} or Subsection 3.3 in \cite{kp}). To address this challenge, we provide an exponential estimate that enables us to establish the required exponential tightness. In comparison with the weak convergence approach used in \cite{q1, q2}, our method does not require the strong averaging principle, resulting in a more concise proof of the LDP.

The contributions of this paper are twofold. First, we establish the LDP for the slow component at any fixed time in the context of fully coupled multiscale multivalued stochastic systems. Moreover, our derived rate function coincides with that of Theorem 3.7 in \cite{q1}. The LDP for the slow component in path space cannot be proved in this case. This limitation arises because, for the multiscale system (\ref{Eq1i}), it is not possible to prove the existence of the limit
\ce
\lim\limits_{\varepsilon\rightarrow 0}\varepsilon\log\mathbb{E}\left[\exp\left(-\frac{1}{\varepsilon}\left(h_1(X^{\varepsilon,\gamma}_{t_1})+\cdots+h_k(X^{\varepsilon,\gamma}_{t_k})\right)\right)\right]
\de
for any $0\leq t_1\leq \cdots \leq t_k$ and $h_1, \cdots, h_k\in C_b( \overline{\mathcal{D}(A)})$, which is a crucial requirement for the family $\{X^{\varepsilon,\gamma}, \varepsilon, \gamma > 0\}$ to satisfy the LDP in $C([0,T], \overline{\mathcal{D}(A)})$ (see, e.g., \cite[Corollary 4.29]{fk}). Second, we prove a comparison principle for first-order Hamilton-Jacobi-Bellman equations with multivalued maximal monotone operators under weaker conditions than those assumed in existing works (cf. \cite{hI, kS}).

The remainder of this paper is organized as follows. In the next section, we introduce the necessary notation, background on maximal monotone operators, large deviation principles, and Hamilton-Jacobi-Bellman equations. Section \ref{main} presents the main result. The proof is provided in Section \ref{ldptimeproo}. An illustrative example is given in Section \ref{exam}. Finally, Section \ref{app} collects auxiliary results, including properties of $\bar{b}_1$ and $\bar{a}_1$, as well as a comparison principle for a general Hamilton-Jacobi-Bellman equation.

Throughout the paper, the following convention is adopted: $C$, with or without subscripts, denotes positive constants whose values may vary from one occurrence to another.

\section{Preliminary}

In this section, we introduce notations, maximal monotone operators, large deviation principles and Hamilton-Jacobi-Bellman equations.

\subsection{Notations}\label{nn}

Let $|\cdot|, \|\cdot\|$ be the norms of a vector and a matrix, respectively. Let $\langle\cdot,\cdot\rangle$ be the inner product of vectors on $\mR^n$. $U^{*}$ denotes the transpose of the matrix $U$.

Let $C(\mR^n)$ be the set of all  functions which are continuous on $\mR^n$. Let $C_b(\mR^n)$ be the set of all  continuous and bounded functions. Let ${\rm Lip}_b(\mR^n)$ be the set of all Lipschitz continuous and bounded functions. The function space $C^{2}(\mR^n)$ represents the collection of all functions in $C(\mR^n)$ with continuous derivatives of order up to 2. The function space $C_b^2(\mR^n)$ stands for the subspace of $C^2(\mR^n)$, consisting of functions whose derivatives up to order 2 and themselves are bounded. The function space $C_c^2(\mR^n)$ is the subspace of $C^2(\mR^n)$, consisting of functions with compact support.

\subsection{Maximal monotone operators}\label{mmo}

For a multivalued operator $A: \mR^n\mapsto 2^{\mR^n}$, where $2^{\mR^n}$ stands for all the subsets of $\mR^n$, set
\ce
&&\cD(A):= \left\{x\in \mR^n: A(x) \ne \emptyset\right\},\\
&&Gr(A):= \left\{(x,y)\in \mR^{2n}:x \in \cD(A), ~ y\in A(x)\right\}.
\de
We say that $A$ is monotone if $\langle x_1 - x_2, y_1 - y_2 \rangle \geq 0$ for any $(x_1,y_1), (x_2,y_2) \in Gr(A)$, and $A$ is maximal monotone if 
$$
(x_1,y_1) \in Gr(A) \iff \langle x_1-x_2, y_1 -y_2 \rangle \geq 0, \quad \forall (x_2,y_2) \in Gr(A).
$$
For example, for a lower semicontinuous convex function $\psi:\mR^n\mapsto(-\infty, +\infty]$ with ${\rm Int}(Dom(\psi))\neq \emptyset$, where $Dom(\psi)\equiv\{x\in\mR^n; \psi(x)<\infty\}$ and $\operatorname{Int}(Dom(\psi))$ is the interior of $Dom(\psi)$, we define the subdifferential operator of the function $\psi$:
$$
\partial\psi(x):=\{y\in\mR^n: \<y,z-x\>+\psi(x)\leq \psi(z), \forall z\in\mR^n\}.
$$
Then $\partial\psi$ is a maximal monotone operator. 

In the following, we recall some properties of a maximal monotone operator $A$ (cf.\cite{cepa1}):
\begin{enumerate}[(i)]
\item
${\rm Int}(\cD(A))$ and $\overline{\mathrm{\cD}(A)}$ are convex subsets of $\mR^n$ with ${\rm Int}\left( \overline{\mathrm{\cD}(A)} \right) = {\rm Int}\( \mathrm{\cD}(A) \) 
$, where ${\rm Int}(\cD(A))$ denotes the interior of the set $\cD(A)$. 
\item $A$ is locally bounded on ${\rm Int}(\cD(A))$, that is, for any compact subset $\Gamma$ of ${\rm Int}(\cD(A))$, $\underset{x\in \Gamma}{\cup}A(x)$ is bounded. 
\end{enumerate}

Take any $T>0$ and fix it. Let $\sV_{0}$ be the set of all continuous functions $K: [0,T]\mapsto\mR^n$ with finite variations and $K_{0} = 0$. For $K\in\sV_0$ and $s\in [0,T]$, we shall use $|K|_{0}^{s}$ to denote the variation of $K$ on $[0,s]$ and write $|K|_{TV}:=|K|_{0}^{T}$. Set
\ce
&&\sA:=\Big\{(X,K): X\in C([0,T],\overline{\cD(A)}), K \in \sV_0, \\
&&\qquad\qquad\quad~\mbox{and}~\langle X_{t}-x, \dif K_{t}-y\dif t\rangle \geq 0 ~\mbox{for any}~ (x,y)\in Gr(A)\Big\}.
\de
And about $\sA$ we have the following two results (cf.\cite[Proposition 4.1 and 4.4]{cepa2} or \cite[Proposition 3.3 and 3.4]{ZXCH}).

\bl\label{equi}
For $X\in C([0,T],\overline{\cD(A)})$ and $K\in \sV_{0}$, the following statements are equivalent:
\begin{enumerate}[(i)]
	\item $(X,K)\in \sA$.
	\item For any $(x,y)\in C([0,T],\mR^{2n})$ with $(x_t, y_t)\in Gr(A)$, it holds that 
	$$
	\left\langle X_t-x_t, \dif K_t-y_t\dif t\right\rangle \geq0.
	$$
	\item For any $(X^{'},K^{'})\in \sA$, it holds that 
	$$
	\left\langle X_t-X_t^{'},\dif K_t-\dif K_t^{'}\right\rangle \geq0.
	$$
\end{enumerate}
\el

\bl\label{inteineq}
Assume that $\text{Int}(\cD(A))\ne\emptyset$. For any $a\in \text{Int}(\cD(A))$, there exist $M_1 >0$, and $M_{2},M_{3}\geq0$ such that  for any $(X,K)\in \sA$ and $0\leq s<t\leq T$,
$$
\int_s^t{\left< X_r-a, \dif K_r \right>}\geq M_1\left| K \right|_{s}^{t}-M_2\int_s^t{\left| X_r-a\right|}\dif r-M_3\left( t-s \right) .
$$
\el

\subsection{Large deviation principles}\label{ldp} (\cite[Chapter 3 and 4]{fk})

Let $(\mS,\rho)$ be a Polish space. For each $\e>0$, let $X^{\e}$ be a $\mS$-valued random variable given on $(\Omega, \mathscr{F}, \{\mathscr{F}_t\}_{t\in[0,T]}, \mP)$.

\bd\label{compleve}
The function $I: \mathbb{S}\mapsto [0,\infty]$ is called a rate function if $I$ is lower semicontinuous. Moreover, a rate function $I$ is called a good rate function if for each $M\in\mR$, $\{\varsigma\in \mathbb{S}:I(\varsigma)\leq M\}$ is a compact subset of $\mathbb{S}$.
\ed

\bd
We say that $\{X^{\e}\}$ satisfies the large deviation principle with the speed $\e^{-1}$ and the good rate function $I$, if for any subset $B\in \sB(\mS)$,
$$
-\inf\limits_{\varsigma\in {\rm Int}(B)}I(\varsigma)\leq\liminf_{\e\rightarrow 0}\e\log\mP(X^{\e}\in{\rm Int}(B))\leq \limsup\limits_{\e\rightarrow 0}\e\log\mP(X^{\e}\in \bar{B})\leq -\inf\limits_{\varsigma\in \bar{B}}I(\varsigma),
$$
where ${\rm Int}(B)$ and $\bar{B}$ denote the interior and the closure of $B$, respectively and they are taken in $\mS$.
\ed

\bd
$\{X^{\e}\}$ is said to be exponentially tight if for every $r>0$, there exists a compact subset $\Gamma_r$ of $\mS$ such that
$$
\limsup\limits _{\e\rightarrow 0}\e \log \mP\left(X^\e\notin \Gamma_r\right) \leq -r.
$$
\ed

\subsection{Hamilton-Jacobi-Bellman equations}\label{hjbeqconv}

Let $\Lambda$ be an index set and
\ce
&&H_{\varepsilon}(z, p, Q): (\overline{\cD(A)}\times\mR^m)\times \mathbb{R}^{(m+n)} \times \mM_{(m+n)}\mapsto \mathbb{R},\\
&&H_i(x, p; \l): \overline{\cD(A)} \times \mathbb{R}^n\times \Lambda \mapsto \mathbb{R}, \quad i=0,1
\de
be continuous, where $\mM_{(m+n)}$ denotes the collection of real $m+n$ order symmetric matrices and $z:=(x,y)$ for $x\in\overline{\cD(A)}$ and $y\in\mR^m$. We also define some following domains
$$
\begin{aligned}
& D_{\varepsilon,+}:=\{f: f \in C^2(\overline{\cD(A)}\times\mR^m), f \text { has compact finite level sets }\},\\
& D_{\varepsilon,-}:=-D_{\varepsilon,+},\\
&D_+:=\{f: f \in C^2(\overline{\cD(A)}), f \text { has compact finite level sets }\},\\
&D_-:=-D_+.
\end{aligned}
$$

In the following, we consider the following Hamilton-Jacobi-Bellman equations:
\be
&&\partial_t u_{\varepsilon}(t, z)\in H_{\varepsilon}\left(z, \p_z u_{\varepsilon}(t, z), \p_{zz}^2 u_{\varepsilon}(t, z)\right)-\<A(x),\p_x u_\e(t,z)\>, \label{hehjb}\\
&&\partial_t u(t, x) \in \inf _{\l \in \Lambda} H_0\left(x, \p_x u(t, x); \l\right)-\<A(x),\p_x u(t,x)\>,\label{h0hjb}\\
&&\partial_t u(t, x) \in \sup _{\l \in \Lambda} H_1\left(x, \p_x u(t, x); \l\right)-\<A(x),\p_x u(t,x)\>,\label{h1hjb}
\ee
where $\p_z u_{\varepsilon}(t, z)$ and $\p_{zz}^2 u_{\varepsilon}(t, z)$ denote the first and second partial derivatives of $u_{\varepsilon}(t, z)$ in $z$. So, we give the definitions of viscosity subsolutions, supersolutions and solutions for Eq.(\ref{hehjb}), Eq.(\ref{h0hjb}) and Eq.(\ref{h1hjb}) (\cite{ffk, Za}).

\bd\label{viscsoludefi} 
$(i)$ We call a bounded measurable function $u$ a viscosity subsolution to Eq.(\ref{h0hjb}) if $u$ is upper semicontinuous, and for each
$$
u_0(t, x)=\phi(t)+f_0(x), \quad \phi \in C^1([0,T]), f_0 \in D_{+},
$$
and each $(t_0,x_0)\in [0,T]\times\overline{\cD(A)}$ satisfying $u-u_0$ has a local maximum at $(t_0,x_0)$, 
\ce
\partial_t u_0\left(t_0, x_0\right)\leq\inf _{\l \in \Lambda} H_0\left(x_0, \p_x u_0(t_0, x_0); \l\right)-A_{\ast}(x_0,\p_x u_0\left(t_0, x_0\right)),
\de
where for $x\in\overline{\cD(A)}$ and $v\in\mR^n$
\ce
A_{\ast}(x, v):=\liminf\limits_{(x', w) \rightarrow (x,v), \zeta\in A(x')}\langle \zeta, w\rangle.
\de

$(ii)$ We call a bounded measurable function $u$ a viscosity supersolution to Eq.(\ref{h1hjb}) if $u$ is lower semicontinuous, and for each
$$
u_1(t, x)=\phi(t)+f_1(x), \quad \phi \in C^1([0,T]), f_1 \in D_{-},
$$
and each $(t_0, x_0) \in [0,T]\times\overline{\cD(A)}$ satisfying $u-u_1$ has a local minimum at $(t_0,x_0)$, 
\ce
\partial_t u_1\left(t_0, x_0\right)\geq\sup _{\l \in \Lambda} H_1\left(x_0, \p_x u_1(t_0, x_0); \l\right)-A^{\ast}(x_0,\p_x u_1\left(t_0, x_0\right)),
\de
where for $x\in\overline{\cD(A)}$ and $v\in\mR^n$
\ce
A^{\ast}(x, v):=\limsup\limits_{(x', w) \rightarrow (x,v), \zeta\in A(x')}\langle \zeta, w\rangle.
\de
\ed

By minor modification, we can define viscosity subsolutions and supersolutions for Eq.(\ref{hehjb}). If a function is both a viscosity subsolution as well as a viscosity supersolution for Eq.(\ref{hehjb}), it is a viscosity solution.

\br\label{Aastsemicont}
$A_*$ is lower semicontinuous and $A^*$ is upper semicontinuous (cf. \cite{Za}).
\er

Besides, we define a class of compact sets in $\overline{\cD(A)}\times\mR^m$ by
$$
\cQ:=\{K\times\tilde{K}: K\subset\subset\overline{\cD(A)}, \tilde{K}\subset\subset\mR^m\},
$$
where $K\subset\subset\overline{\cD(A)}$ means that $K$ is a compact subset of $\overline{\cD(A)}$. Set
\ce
H_\varepsilon f(z):=H_{\varepsilon}\left(z, \p_z f(z), \p_{zz}^2 f(z)\right), \quad f\in D_{\e,+}\cup D_{\e,-}.
\de
Then we give two conditions.

\bco\label{condsup} 
For each $f_0 \in D_{+}$and $\l \in \Lambda$, there exists $f_{0, \varepsilon} \in D_{\e,+}$ (which may depend on $\l$) such that

$(i)$ for each $c>0$, there exists $K \times \tilde{K} \in \mathcal{Q}$ satisfying
$$
\left\{(x, y): H_\e f_{0, \e}(x, y) \geq-c\right\} \cap\left\{(x, y): f_{0, \e}(x, y) \leq c\right\} \subset K \times \tilde{K};
$$

$(ii)$ for each $K \times \tilde{K} \in \mathcal{Q}$,
$$
\lim _{\e \rightarrow 0} \sup _{(x, y) \in K \times \tilde{K}}\left[\left|f_{0, \e}(x, y)-f_0(x)\right|+\left|\p_xf_{0, \e}(x, y)-\p_xf_0(x)\right|\right]=0;
$$

$(iii)$ whenever $\left(x_\e, y_\e\right) \in K \times \tilde{K} \in \mathcal{Q}$ satisfies $x_\e \rightarrow x$,
$$
\limsup _{\e \rightarrow 0} H_\e f_{0, \e}\left(x_\e, y_\e\right) \leq H_0\left(x, \p_x f_0(x); \l\right).
$$
\eco

\bco\label{condinf}
For each $f_1 \in D_{-}$and $\l \in \Lambda$, there exists $f_{1, \e} \in D_{\e,-}$ (which may depend on $\l$) such that

$(i)$ for each $c>0$, there exists $K \times \tilde{K} \in \mathcal{Q}$ satisfying
$$
\left\{(x, y): H_\e f_{1, \e}(x, y) \leq c\right\} \cap\left\{(x, y): f_{1, \e}(x, y) \geq-c\right\} \subset K \times \tilde{K};
$$

$(ii)$ for each $K \times \tilde{K} \in \mathcal{Q}$,
$$
\lim _{\e \rightarrow 0} \sup _{(x, y) \in K \times \tilde{K}}\left[\left|f_{1, \e}(x, y)-f_1(x)\right|+\left|\p_xf_{1, \e}(x, y)-\p_xf_1(x)\right|\right]=0;
$$

$(iii)$ whenever $\left(x_\e, y_\e\right) \in K \times \tilde{K} \in \mathcal{Q}$, and $x_\e \rightarrow x$,
$$
\liminf _{\e \rightarrow 0} H_\e f_{1, \e}\left(x_\e, y_\e\right) \geq H_1\left(x, \p_x f_1(x); \l\right).
$$
\eco

Next, let $u_{\e}$ be a viscosity solution to Eq.(\ref{hehjb}). We define the limits of $u_{\e}$ as $\e\rightarrow 0$.

\bd\label{ueglimidefi}
\ce
u_{\uparrow}(t,x):=&\sup\Big\{\limsup\limits_{\e\rightarrow 0}u_{\e}(t_\e,x_\e,y_\e): \exists (t_\e,x_\e,y_\e)\in[0,T]\times K\times\tilde{K}, \\
&(t_\e,x_\e)\rightarrow (t,x), K\times\tilde{K}\in\cQ\Big\},\\
u_{\downarrow}(t,x):=&\inf\Big\{\liminf\limits_{\e\rightarrow 0}u_{\e}(t_\e,x_\e,y_\e): \exists (t_\e,x_\e,y_\e)\in[0,T]\times K\times\tilde{K}, \\
&(t_\e,x_\e)\rightarrow (t,x), K\times\tilde{K}\in\cQ\Big\}.
\de
Define $\bar{u}$ to be the upper semicontinuous regularization of $u_{\uparrow}$, and $\underline{u}$ the lower semicontinuous regularization of $u_{\downarrow}$.
\ed

Finally, we have the following two results.

\bl\label{supersubsolu}
Suppose that the viscosity solutions $u_\e$ to Eq.(\ref{hehjb}) are uniformly bounded. Then under Condition \ref{condsup}, $\bar{u}$ is a viscosity subsolution of Eq.(\ref{h0hjb})
and under Condition \ref{condinf}, $\underline{u}$ is a viscosity supersolution of Eq.(\ref{h1hjb}).
\el
\begin{proof}
Since the proofs of a viscosity subsolution for Eq.(\ref{h0hjb}) and a viscosity supersolution for Eq.(\ref{h1hjb}) are similar, we only prove that under Condition \ref{condsup}, $\bar{u}$ is a viscosity subsolution of Eq.(\ref{h0hjb}).

Take
$$
u_0(t, x)=\phi(t)+f_0(x), \quad \phi \in C^1([0,T]), f_0 \in D_{+},
$$
and $(t_0,x_0)\in [0,T]\times\overline{\cD(A)}$ satisfying $\bar{u}-u_0$ has a local maximum at $(t_0,x_0)$. Then we change $\phi$ and $f_0$ such that $(t_0,x_0)$ is a unique global maximum (cf. \cite{ffk}).   Let $\l \in \Lambda$ be given. So, by Condition \ref{condsup}, there exists a $f_{0,\e}\in D_{\e,+}$ such that $(i)$-$(iii)$ in Condition \ref{condsup} hold. Let $u_{0,\e}(t,z):=\phi(t)+f_{0,\e}(z)$. Since $u_{\varepsilon}(t,z)$ is bounded, and $u_{0, \varepsilon}(t,z)$ has compact level sets, there exists $\left(t_{\varepsilon}, z_{\varepsilon}\right)=(t_\e,x_\e,y_\e) \in[0, T] \times \overline{\cD(A)}\times\mR^m$ such that $u_{\varepsilon}-u_{0, \varepsilon}$ has a local maximum at $(t_\e,z_\e)$. Note that $u_{\varepsilon}$ is a viscosity subsolution of Eq.(\ref{hehjb}). Thus, it holds that
\ce
\partial_t u_{0,\varepsilon}(t_\e, z_\e)\leq H_{\varepsilon}\left(z_\e, \p_z u_{0,\varepsilon}(t_\e, z_\e), \p_{zz}^2 u_{0,\varepsilon}(t_\e, z_\e)\right)-A_{\ast}(x_\e, \p_x u_{0,\e}(t_\e,z_\e)),
\de
and 
\be
\partial_t \phi(t_\e)\leq H_\e f_{0,\e}(z_\e)-A_{\ast}(x_\e, \p_x f_{0,\e}(z_\e)).
\label{phief0e}
\ee
From the above deduction, we infer that for any $c>0$, $H_{\varepsilon} f_{0, \varepsilon}\left(z_{\varepsilon}\right)\geq -c$ and $f_{0, \varepsilon}\left(z_{\varepsilon}\right)\leq c$. So, by $(i)$ in Condition \ref{condsup}, there exists $K \times \tilde{K} \in \mathcal{Q}$ such that $z_{\varepsilon}=\left(x_{\varepsilon}, y_{\varepsilon}\right) \in K \times \tilde{K}$.

Since $K \times \tilde{K}$ is compact in $\overline{\cD(A)}\times\mR^m$, there exists a subsequence denoted still by $\left\{\left(t_{\varepsilon}, z_{\varepsilon}\right)\right\}$ and a $\left(\hat t_0, \hat x_0\right) \in[0, T] \times \overline{\cD(A)}$ such that $t_{\varepsilon} \rightarrow \hat{t}_0$ and $x_{\varepsilon} \rightarrow \hat{x}_0$. We claim that $\left(\hat t_0, \hat x_0\right)=(t_0,x_0)$. Indeed, since $(t_\e,x_\e)\rightarrow \left(\hat t_0, \hat x_0\right)$ and $u_{\varepsilon}-u_{0, \varepsilon}$ has a local maximum at $(t_\e,z_\e)$, $\bar{u}-u_0$ has a local maximum at $\left(\hat t_0, \hat x_0\right)$ and $\left(\hat t_0, \hat x_0\right)=(t_0,x_0)$.

Finally, taking the superior limit on two sides of (\ref{phief0e}), by $(ii)$-$(iii)$ in Condition \ref{condsup} and Remark \ref{Aastsemicont}, we obtain that
\ce
\partial_t \phi(t_0)\leq H_0\left(x_0, \p_x f_0\left(x_0\right); \l\right)-A_{\ast}(x_0, \p_x f_{0}(x_0)),
\de 
that is,
\ce
\p_t u_0(t_0,x_0)\leq H_0\left(x, \p_x u_0(t_0, x_0); \l\right)-A_{\ast}(x_0,\p_x u_0(t_0,x_0)\>.
\de
Note that $t_0, x_0$ and $u_0$ are all chosen independent of $\l$. Thus, we take the infimum with respect to $\l$ on both sides and conclude that
$$
\p_t u_0(t_0,x_0)\leq \inf _{\l \in \Lambda}H_0\left(x, \p_x u_0(t_0, x_0); \l\right)-A_{\ast}(x_0,\p_x u_0(t_0,x_0)\>.
$$
The proof is complete.
\end{proof}

\bl\label{uehlimi}
Suppose that the assumptions in Lemma \ref{supersubsolu} hold. If a comparison principle between the viscosity subsolutions to Eq.(\ref{h0hjb}) and the viscosity supersolutions to Eq.(\ref{h1hjb}) holds, then $u:=\bar{u}=\underline{u}$ and 
$$
\lim\limits_{\e\rightarrow 0}\sup\limits_{t\in[0,T]}\sup\limits_{(x,y)\in K \times \tilde{K}}|u_\e(t,x,y)-u(t,x)|=0, \quad \forall K \times \tilde{K}\in\mathcal{Q}.
$$
\el
\begin{proof}
By Lemma \ref{supersubsolu}, we know that $\bar{u}$ is a viscosity subsolutions to Eq.(\ref{h0hjb}) and $\underline{u}$ is a viscosity supersolutions to Eq.(\ref{h1hjb}). Since a comparison principle between the viscosity subsolutions to Eq.(\ref{h0hjb}) and the viscosity supersolutions to Eq.(\ref{h1hjb}) holds, it follows that
$$
\bar{u}\leq \underline{u}.
$$
Besides, by Definition \ref{ueglimidefi}, it is obvious that 
$$
\bar{u}\geq \underline{u}.
$$
Collecting the above deduction, we conclude that $u:=\bar{u}=\underline{u}$, which together with Definition \ref{ueglimidefi} yields that
$$
\lim\limits_{\e\rightarrow 0}\sup\limits_{t\in[0,T]}\sup\limits_{(x,y)\in K \times \tilde{K}}|u_\e(t,x,y)-u(t,x)|=0, \quad \forall K \times \tilde{K}\in\mathcal{Q}.
$$
The proof is complete.
\end{proof}

\section{Main results}\label{main}

In this section, we formulate the main result in this paper.

We recall the system (\ref{Eq1i}), i.e.
\ce\left\{\begin{array}{l}
\dif X_{t}^{\e,\g}\in -A(X_{t}^{\e,\g})\dif t+b_{1}(X_{t}^{\e,\g},Y_{t}^{\e,\g})\dif t+\sqrt{\e}\s_{1}(X_{t}^{\e,\g},Y_{t}^{\e,\g})\dif W^1_{t},\\
X_{0}^{\e,\g}=x_0\in\overline{\cD(A)},\quad  0\leq t\leq T,\\
\dif Y_{t}^{\e,\g}=\frac{1}{\g}b_{2}(X_{t}^{\e,\g},Y_{t}^{\e,\g})\dif t+\frac{1}{\sqrt{\g}}\s_{2}(X_{t}^{\e,\g},Y_{t}^{\e,\g})\dif W^2_{t},\\
Y_{0}^{\e,\g}=y_0,\quad  0\leq t\leq T.
\end{array}
\right.
\de

Assume:
\begin{enumerate}[$(\mathbf{H}_{A})$]
\item $0\in{\rm Int}(\cD(A))$.
\end{enumerate}
\begin{enumerate}[$(\mathbf{H}^1_{b_{1}, \s_{1}})$]
\item
There exists a constant $L_{b_{1}, \s_{1}}>0$ such that for $x_{i}\in\overline{\cD(A)}$, $y_{i}\in\mR^m$, $i=1, 2$,
\ce
|b_{1}(x_{1},y_{1})-b_{1}(x_{2},y_{2})|^{2}+\|\s_{1}(x_1,y_{1})-\s_{1}(x_2,y_{2})\|^{2}\leq L_{b_{1},\s_{1}}\(|x_{1}-x_{2}|^{2}+|y_{1}-y_{2}|^{2}\),
\de
and for $x\in\overline{\cD(A)}, y\in\mR^m$
\ce
|b_1(x,y)|^2\leq L_{b_{1},\s_{1}},\quad  \|\s_{1}(x,y)\|^2\leq L_{b_{1},\s_{1}}.
\de
\end{enumerate}
\begin{enumerate}[$(\mathbf{H}^2_{b_{1}, \s_{1}})$]
\item
$\partial_y b_1(x,y), \partial_{yy} b_1(x,y), \partial_x\s_1\s^*_1(x,y), \partial_y\s_1\s^*_1(x,y), \partial_{yy}\s_1\s^*_1(x,y) $ exist for any $(x,y)\in\overline{\cD(A)}\times\mR^m$, and $\partial_x\s_1\s^*_1(x,y), \partial_y\s_1\s^*_1(x,y)$ are uniformly bounded.
\end{enumerate}
\begin{enumerate}[$(\mathbf{H}^1_{b_{2}, \s_{2}})$]
\item
There exists a constant $L_{b_{2}, \s_{2}}>0$ such that for $x_{i}\in\overline{\cD(A)}$, $y_{i}\in\mR^m$, $i=1, 2$,
\ce
|b_{2}(x_{1},y_{1})-b_{2}(x_{2},y_{2})|^{2}+\|\s_{2}(x_{1},y_{1})-\s_{2}(x_{2},y_{2})\|^{2}\leq L_{b_{2}, \s_{2}}\(|x_{1}-x_{2}|^{2}+|y_{1}-y_{2}|^{2}\).
\de
\end{enumerate}
\begin{enumerate}[$(\mathbf{H}^2_{b_{2}, \s_{2}})$]
\item
There exists a constant $\b>0$ satisfying $\b>2L_{b_2,\sigma_2}$ such that for $x\in\overline{\cD(A)}$, $y_{i}\in\mR^m$, $i=1, 2$,
\ce
2\<y_{1}-y_{2},b_{2}(x,y_{1})-b_{2}(x,y_{2})\>
+\|\s_{2}(x,y_{1})-\s_{2}(x,y_{2})\|^{2}\leq -\b|y_{1}-y_{2}|^{2}.
\de
\end{enumerate}
\begin{enumerate}[$(\mathbf{H}^3_{b_2,\s_2})$]
\item There exists a positive function $\zeta(\cdot)\in C^2(\mR^m)$ such that $\zeta$ has compact finite level sets and there exist two constants $L_1, L_2>0$ and a compact set $B\subset\mR^m$ such that for any $x\in\overline{\cD(A)}$, 
\ce
\sL^x_1 \zeta(y)\leq -L_1\zeta(y)+L_2 I_B(y),
\de
where
\ce
(\sL^x_1\zeta)(y):=\<b_2(x,y),\p_y\zeta(y)\>+\frac{1}{2}{\rm tr}\left(\s_2\s^*_2(x,y)\p^2_{yy}\zeta(y)\right).
\de
\end{enumerate}
\begin{enumerate}[$(\mathbf{H}^4_{b_{2}, \s_{2}})$]
\item
$\partial_y b_2(x,y), \partial_{yy} b_2(x,y), \partial_y\s_2(x,y), \partial_{yy}\s_2(x,y)$ exist for any $(x,y)\in\overline{\cD(A)}\times\mR^m$.
\end{enumerate}

\br
$(i)$ $(\mathbf{H}^1_{b_{2}, \s_{2}})$ implies that there exists a constant $\bar{L}_{b_{2}, \s_{2}}>0$ such that for $x\in\overline{\cD(A)}$, $y\in\mR^m$,
\be
|b_{2}(x,y)|^{2}+\|\s_{2}(x,y)\|^{2}\leq \bar{L}_{b_{2}, \s_{2}}(1+|x|^{2}+|y|^{2}).
\label{b2nu}
\ee

$(ii)$ $(\mathbf{H}^1_{b_{2}, \s_{2}})$ and $(\mathbf{H}^2_{b_{2}, \s_{2}})$ yield that there exists a constant $C>0$ such that for $x\in\overline{\cD(A)}$, $y\in\mR^m$
\be
&&2\<y,b_{2}(x,y)\>+\|\s_{2}(x,y)\|^{2}\leq -\a|y|^{2}+C(1+|x|^{2}),
\label{bemu}
\ee
where $\a:=\b-2L_{b_{2}, \s_{2}}>0$.

$(iii)$ $(\mathbf{H}_{A})$, $(\mathbf{H}^1_{b_{1}, \s_{1}})$ and $(\mathbf{H}^1_{b_{2}, \s_{2}})$ assure the well-posedness of the system (\ref{Eq1i}). By $(\mathbf{H}^1_{b_{2}, \s_{2}})$  and $(\mathbf{H}^2_{b_{2}, \s_{2}})$, we get the existence and uniqueness of the invariant measures for the frozen equation. In terms of $(\mathbf{H}_{A})$, $(\mathbf{H}^1_{b_{1}, \s_{1}})$-$(\mathbf{H}^2_{b_{1}, \s_{1}})$ and $(\mathbf{H}^1_{b_{2}, \s_{2}})$-$(\mathbf{H}^4_{b_{2}, \s_{2}})$, the large deviation principle is established. 
\er

Under $(\mathbf{H}_{A})$, $(\mathbf{H}^1_{b_{1}, \s_{1}})$ and $(\mathbf{H}^1_{b_{2}, \s_{2}})$, by \cite[Theorem 3.2]{q2} we know that the system (\ref{Eq1i}) has the unique strong solution $(X^{\e,\g}, K^{1,\e,\g}, Y^{\e,\g})$. That is, $(X^{\e,\g}, K^{1,\e,\g})\in\sA$ and 
\ce\left\{\begin{array}{l}
X_{t}^{\e,\g}=x_0-K^{1,\e,\g}_t+\int_0^t b_{1}(X_{s}^{\e,\g},Y_{s}^{\e,\g})\dif s+\sqrt{\e}\int_0^t\s_{1}(X^{\e,\g},Y_{s}^{\e,\g})\dif W^1_{s},\\
Y_{t}^{\e,\g}=y_0+\frac{1}{\g}\int_0^t b_{2}(X_{s}^{\e,\g},Y_{s}^{\e,\g})\dif s+\frac{1}{\sqrt{\g}}\int_0^t\s_{2}(X_{s}^{\e,\g},Y_{s}^{\e,\g})\dif W^2_{s}.
\end{array}
\right.
\de

Fix $x\in\overline{\cD(A)}$, and consider the following SDE:
\be\left\{\begin{array}{l}
\dif Y_{t}^{x}=b_{2}(x,Y_{t}^{x})\dif t+\s_{2}(x,Y_{t}^{x})\dif W^2_{t},\\
Y_{0}^{x}=y_0,\quad  t\geq 0.
\end{array}
\right.
\label{frozequa}
\ee
Under $(\mathbf{H}^1_{b_{2}, \s_{2}})$, the above equation has the unique strong solution $Y^{x,y_0}$ and $Y^{x,y_0}$ is a Markov process with the transition semigroup $\{P_{t}^{x}\}_{t\geq0}$ and the infinitesimal generator $\sL^x_1$. In addition, under $(\mathbf{H}^1_{b_{2}, \s_{2}})$ and $(\mathbf{H}^2_{b_{2}, \s_{2}})$, $Y^{x,y_0}$ has the unique invariant probability measure $\mu^x$ (\cite{WangFY}). 

Next, we define
\ce
\bar{b}_1(x):=\int_{\mR^m}b_1(x,y)\mu^x(\dif y), \quad \bar{a}_1(x):=\int_{\mR^m}\s_1\s^*_1(x,y)\mu^x(\dif y).
\de
Consider the following Hamilton-Jacobi-Bellman equation
\be\left\{\begin{array}{l} 
\p_t u(t,x)\in \bar{H}(x,\p_x u(t,x)) -\<A(x),\p_x u(t,x)\>,\\
u(0,x)=h(x), \quad x\in\overline{\cD(A)},
\end{array}
\right.
\label{hjblimi}
\ee
where 
\ce
\bar{H}(x,p):=\<\bar{b}_1(x),p\>-\frac{1}{2}\<\bar{a}_1(x)p,p\>, \quad x\in\overline{\cD(A)},\quad  p\in\mR^n,
\de
and $h$ is a continuous bounded function on $\overline{\cD(A)}$. Based on the matrix theory, for the positive semidefinite matrix $\bar{a}_1(x)$, there exists a positive semidefinite matrix $\bar{\s}_1(x)$ such that $\bar{a}_1(x)=\bar{\s}_1(x)\bar{\s}_1(x)$ for $ x\in\overline{\cD(A)}$. Then Eq.(\ref{hjblimi}) is rewritten as
\ce\left\{\begin{array}{l} 
\p_t u(t,x)-\<\bar{b}_1(x),\p_x u(t,x)\>+\frac{1}{2}|\bar{\s}_1(x)\p_x u(t,x)|^2\in -\<A(x),\p_x u(t,x)\>,\\
u(0,x)=h(x), \quad x\in\overline{\cD(A)}.
\end{array}
\right.
\de
Note that
\ce
\frac{1}{2}|\bar{\s}_1(x)p|^2=\sup\limits_{z\in\mR^n}\left\{-\<p,\bar{\s}_1(x)z\>-\frac{|z|^2}{2}\right\}.
\de
Thus, by Theorem 4 in \cite{dT}, we know that 
\be
u_0^h(t,x_0)&:=&\inf\left\{\frac{1}{2}\int_0^t|z(s)|^2\dif s+h(X^z_{x_0}(t)): z\in L^2([0,T],\mR^n)\right\}\label{hjblimisolu}\\
&=&\inf\limits_{x\in\overline{\cD(A)}}\left\{\inf\left\{\frac{1}{2}\int_0^t|z(s)|^2\dif s: z\in L^2([0,T],\mR^n), s.t. X^z_{x_0}(t)=x\right\}+h(x)\right\}\no
\ee
is a viscosity solution of Eq.(\ref{hjblimi}), where $(X^z_{x_0},K^z_{x_0})$ is a solution of the following multivalued differential equation
\be\left\{\begin{array}{l}
\dif X^z_{x_0}(t)\in -A(X^z_{x_0}(t))\dif t+\bar{b}_{1}(X^z_{x_0}(t))\dif t+\bar{\s}_{1}(X^z_{x_0}(t))z(t)\dif t,\quad  0\leq t\leq T,\\
X^z_{x_0}(0)=x_0\in\overline{\cD(A)}.
\end{array}
\right.
\label{ratefuncequa}
\ee

Now, it is the position to state the main result in this paper.

\bt\label{ldptime}
Suppose that $(\mathbf{H}_{A})$, $(\mathbf{H}^1_{b_{1}, \s_{1}})$, $(\mathbf{H}^2_{b_{1}, \s_{1}})$, $(\mathbf{H}^1_{b_{2}, \s_{2}})$-$(\mathbf{H}^4_{b_{2}, \s_{2}})$ hold. If 
\ce
\lim\limits_{\e\rightarrow0}\frac{\g}{\e}=0,
\de
then $\{X^{\e,\g}_t, \e>0\}$ for any $t\in[0,T]$ satisfies the LDP with the speed $\e^{-1}$ and the good rate function $I$ given by
\ce
I(x; x_0, t)=\inf\left\{\frac{1}{2}\int_0^t|z(s)|^2\dif s: z\in L^2([0,T],\mR^n), s.t. X^z_{x_0}(t)=x\right\}.
\de
\et

\br
We mention that the rate function is the same to that in \cite[Theorem 3.7]{q1}.
\er

The proof of Theorem \ref{ldptime} is placed in the next section.

\section{Proof of Theorem \ref{ldptime}}\label{ldptimeproo}

In this section, we follow the line in \cite{ffk} to prove Theorem \ref{ldptime}. That is, we firstly investigate the limit of $-\e\log\mE\left[\exp\left(-\frac{h(X^{\e,\g}_t)}{\e}\right)|X_{0}^{\e,\g}=x, Y_{0}^{\e,\g}=y\right]$ for $h\in C_b(\overline{\cD(A)})$ as $\e\rightarrow0$, show the exponential tightness of $\{X^{\e,\g}_t, \e, \g>0\}$ and by the Bryc formula (\cite[Proposition 3.8]{fk}) establish  the large deviation principle of $\{X^{\e,\g}_t, \e, \g>0\}$ for any $t\in[0,T]$.

\subsection{The Laplace limit at the time $t$}

In this subsection, we observe the limit of $-\e\log\mE\left[\exp\left(-\frac{h(X^{\e,\g}_t)}{\e}\right)|X_{0}^{\e,\g}=x, Y_{0}^{\e,\g}=y\right]$ for $h\in C_b(\overline{\cD(A)})$ as $\e\rightarrow0$. 

To do this, set
\ce
u^h_{\e,\g}(t,x,y):=-\e\log\mE\left[\exp\left(-\frac{h(X^{\e,\g}_t)}{\e}\right)|X_{0}^{\e,\g}=x, Y_{0}^{\e,\g}=y\right], \quad h\in C_b(\overline{\cD(A)}).
\de
For $\Phi\in C_c^2(\overline{\cD(A)}\times\mR^m)$, set 
\ce
(\sL_{\e,\g}\Phi)(x,y):=\<b_1(x,y),\p_x\Phi(x,y)\>+\frac{\e}{2}{\rm tr}\left(\s_1\s^*_1(x,y)\p^2_{xx}\Phi(x,y)\right)+\frac{1}{\g}(\sL^x_1\Phi)(x,y).
\de
So, by the similar deduction to that in \cite[Theorem 4]{Za}, we infer that $u^h_{\e,\g}$ is a unique viscosity solution of the following Cauchy problem:
\be\left\{\begin{array}{ll}
\p_tu\in H_{\e,\g}u-\<A(x),\p_x u\>, \quad\qquad (t,x,y)\in[0,T]\times\overline{\cD(A)}\times\mR^m,\\
u(0,x,y)=h(x), \quad\qquad\quad\quad\qquad (x,y)\in\overline{\cD(A)}\times\mR^m,
\end{array}
\right.
\label{caucprob}
\ee
where 
\ce
&&(H_{\e,\g}u)(t,x,y):=-\e e^{u/\e}\sL_{\e,\g}e^{-u/\e}(t,x,y)\\
&=&\<b_1(x,y),\p_x u(t,x,y)\>-\frac{1}{2}\<\s_1\s^*_1(x,y)\p_x u(t,x,y),\p_x u(t,x,y)\>\\
&&+\frac{\e}{2}{\rm tr}\left(\s_1\s^*_1(x,y)\p^2_{xx}u(t,x,y)\right)+\frac{1}{\g}\[\<b_2(x,y),\p_y u(t,x,y)\>\\
&&\quad +\frac{1}{2}{\rm tr}\left(\s_2\s^*_2(x,y)\p^2_{yy}u(t,x,y)\right)\]-\frac{1}{2\g\e}\<\s_2\s^*_2(x,y)\p_y u(t,x,y),\p_y u(t,x,y)\>.
\de

Next, we study the convergence of the Hamilton-Jacobi-Bellman equation (\ref{caucprob}) in terms of the procedure in Subsection \ref{hjbeqconv}. Set
$$
\begin{aligned}
& D_{\e,+}:=\left\{f; f\in C^2(\overline{\cD(A)}\times\mR^m), f ~\mbox{has a compact finite level set}\right\},\\
& D_{\e,-}:=-D_{\e,+},\\
& D_{+}:=\left\{f\in C^2\left(\overline{\cD(A)}\right): f(x)=\varphi(x)+\log \left(1+|x|^2\right); \varphi \in C_c^2\left(\overline{\cD(A)}\right)\right\}, \\
& D_{-}:=\left\{f\in C^2\left(\overline{\cD(A)}\right): f(x)=\varphi(x)-\log \left(1+|x|^2\right); \varphi \in C_c^2\left(\overline{\cD(A)}\right)\right\}.
\end{aligned}
$$
Taking the index set
$$
\Lambda:=\left\{\l=(\xi, \theta); \xi \in C^2_c(\mR^m), 0<\theta<1\right\}.
$$

{\bf Verification of Condition \ref{condsup}.} For any $f_0\in D_{+}$ and any $\l=(\xi, \t)\in\Lambda$, we take 
\ce
g(y):=\xi(y)+\t\zeta(y),
\de
where $\zeta$ satisfies $(\mathbf{H}^3_{b_2,\s_2})$, and define a perturbed test function
\ce
f_{0,\e}(x,y):=f_0(x)+\g g(y)=f_0(x)+\g\xi(y)+\g\t\zeta(y).
\de
So, 
\ce
(H_{\e,\g}f_{0,\e})(x,y)&=&\<b_1(x,y),\p_x f_0(x)\>-\frac{1}{2}\<\s_1\s^*_1(x,y)\p_x f_0(x),\p_x f_0(x)\>\\
&&+\frac{\e}{2}{\rm tr}\left(\s_1\s^*_1(x,y)\p^2_{xx}f_0(x)\right)+\left[\sL^x_1\xi(y)+\t\sL^x_1\zeta(y)\right]\\
&&-\frac{\g}{2\e}\<\s_2\s^*_2(x,y)\p_y (\xi(y)+\t\zeta(y)),\p_y (\xi(y)+\t\zeta(y))\>.
\de
By the definitions of $D_{+}$ and $(\mathbf{H}^3_{b_2,\s_2})$, we know that $(i)$ and $(ii)$ in Condition \ref{condsup} hold.

Define a function $H_0(x,p;\l): \overline{\cD(A)}\times\mR^n\mapsto\mR$ for any $\l=(\xi, \theta)\in\Lambda$ by
\ce
H_0(x,p;\l):=\sup\limits_{y\in\mR^m}\left\{\<b_1(x,y),p\>-\frac{1}{2}\<\s_1\s^*_1(x,y)p,p\>+\sL^x_1\xi(y)-\t L_1\zeta(y)+\t L_2 I_B(y)\right\}.
\de
Then, whenever $(x_\e, y_\e)\in K\times\tilde{K}\in\cQ$ and $x_\e\rightarrow x$, $(\mathbf{H}^3_{b_2,\s_2})$ and $\lim\limits_{\e\rightarrow0}\g/\e=0$ imply that
\ce
&&\limsup\limits_{\e\rightarrow0}(H_{\e,\g}f_{0,\e})(x_\e,y_\e)\\
&=&\limsup\limits_{\e\rightarrow0}\bigg[\<b_1(x,y_\e),\p_x f_0(x)\>-\frac{1}{2}\<\s_1\s^*_1(x,y_\e)\p_x f_0(x),\p_x f_0(x)\>\\
&&\qquad\qquad+\left[\sL^x_1\xi(y_\e)+\t\sL^x_1\zeta(y_\e)\right]\bigg]\\
&\leq& \sup\limits_{y\in\mR^m}\bigg\{\<b_1(x,y),\p_x f_0(x)\>-\frac{1}{2}\<\s_1\s^*_1(x,y)\p_x f_0(x),\p_x f_0(x)\>+\sL^x_1\xi(y)\\
&&\qquad\quad -\t L_1\zeta(y)+\t L_2 I_B(y)\bigg\}\\
&=& H_0(x,\p_x f_0; \l).
\de
From this, we conclude that $(iii)$ in Condition \ref{condsup} holds.

{\bf Verification of Condition \ref{condinf}.} For any $f_1\in D_{-}$ and any $\l=(\xi, \t)\in\Lambda$, we take 
\ce
&&g(y):=\xi(y)-\t\zeta(y),\\
&&f_{1,\e}(x,y):=f_1(x)+\g g(y)=f_1(x)+\g\xi(y)-\g\t\zeta(y).
\de
Then
\ce
(H_{\e,\g}f_{1,\e})(x,y)&=&\<b_1(x,y),\p_x f_1(x)\>-\frac{1}{2}\<\s_1\s^*_1(x,y)\p_x f_1(x),\p_x f_1(x)\>\\
&&+\frac{\e}{2}{\rm tr}\left(\s_1\s^*_1(x,y)\p^2_{xx}f_1(x)\right)+\left[\sL^x_1\xi(y)-\t\sL^x_1\zeta(y)\right]\\
&&-\frac{\g}{2\e}\<\s_2\s^*_2(x,y)\p_y(\xi(y)-\t\zeta(y)),\p_y(\xi(y)-\t\zeta(y))\>.
\de
Based on the definitions of $D_{-}$ and $(\mathbf{H}^3_{b_2,\s_2})$, we infer that $(i)$ and $(ii)$ in Condition \ref{condinf} hold.

Define a function $H_1(x,p;\l): \overline{\cD(A)}\times\mR^n\mapsto\mR$ for any $\l=(\xi, \theta)\in\Lambda$ by
\ce
H_1(x,p;\l):=\inf\limits_{y\in\mR^m}\left\{\<b_1(x,y),p\>-\frac{1}{2}\<\s_1\s^*_1(x,y)p,p\>+\sL^x_1\xi(y)+\t L_1\zeta(y)-\t L_2 I_B(y)\right\}.
\de
Then, whenever $(x_\e, y_\e)\in K\times\tilde{K}\in\cQ$ and $x_\e\rightarrow x$, $(\mathbf{H}^3_{b_2,\s_2})$ and $\lim\limits_{\e\rightarrow0}\g/\e=0$ imply that
\ce
&&\liminf\limits_{\e\rightarrow0}(H_{\e,\g}f_{1,\e})(x_\e,y_\e)\\
&=&\liminf\limits_{\e\rightarrow0}\bigg[\<b_1(x,y_\e),\p_x f_1(x)\>-\frac{1}{2}\<\s_1\s^*_1(x,y_\e)\p_x f_1(x),\p_x f_1(x)\>\\
&&\qquad\qquad+\left[\sL^x_1\xi(y_\e)-\t\sL^x_1\zeta(y_\e)\right]\bigg]\\
&\geq&\inf\limits_{y\in\mR^m}\bigg\{\<b_1(x,y),\p_x f_1(x)\>-\frac{1}{2}\<\s_1\s^*_1(x,y)\p_x f_1(x),\p_x f_1(x)\>+\sL^x_1\xi(y)\\
&&\qquad\quad+\t L_1\zeta(y)-\t L_2 I_B(y)\bigg\}\\
&=& H_1(x,\p_x f_1; \l).
\de
That is, $(iii)$ in Condition \ref{condsup} is verified.

Combining the above deduction with Lemma \ref{supersubsolu}, we draw the following conclusion.

\bp\label{concsupersubsolu}
Suppose that $(\mathbf{H}_{A})$, $(\mathbf{H}^1_{b_{1}, \s_{1}})$, $(\mathbf{H}^1_{b_{2}, \s_{2}})$ and $(\mathbf{H}^3_{b_{2}, \s_{2}})$ hold. If 
\ce
\lim\limits_{\e\rightarrow0}\frac{\g}{\e}=0,
\de
then $\bar{u}^h$ is a viscosity subsolution of the following Hamilton-Jacobi-Bellman equation
\be\left\{\begin{array}{l}
\p_t u(t,x)\in\inf\limits _{\l \in \Lambda} H_0\left(x, \p_x u(t, x); \l\right)-\<A(x),\p_x u(t,x)\>,\\
u(0,x)=h(x),
\end{array}
\right.
\label{conchohhjb}
\ee
and $\underline{u}^h$ is a viscosity supersolution of the following Hamilton-Jacobi-Bellman equation
\be\left\{\begin{array}{l}
\p_t u(t,x)\in\sup\limits_{\l \in \Lambda}H_1\left(x, \p_x u(t, x); \l\right)-\<A(x),\p_x u(t,x)\>,\\
u(0,x)=h(x),
\end{array}
\right.
\label{conch1hhjb}
\ee
where $\bar{u}^h$ and $\underline{u}^h$ are defined similarly in Definition \ref{ueglimidefi} by replacing $u_{\e}$ with $u^h_{\e,\g}$.
\ep

In the following, for $x\in\overline{\cD(A)}$ and $p\in\mR^n$ set 
$$
H_0(x, p):=\inf\limits _{\l \in \Lambda} H_0(x, p; \l), \quad H_1(x, p):=\sup\limits _{\l \in \Lambda} H_1(x, p; \l).
$$
Then in order to study the relationship among $H_0(x, p), H_1(x, p)$ and $\bar{H}(x,p)$, we first prepare the following result. Consider the following Poisson equation:  for $x\in\overline{\cD(A)}, y\in\mR^m, p\in\mR^n$
\be
(\sL^x_1\kappa)(x,y,p)=-\left[\<b_1(x,y),p\>-\frac{1}{2}\<\s_1\s_1^*(x,y)p, p\>-\<\bar{b}_1(x),p\>+\frac{1}{2}\<\bar{a}_1(x)p, p\>\right].
\label{poisequa}
\ee
Set 
$$
\Psi(x,y,p):=\<b_1(x,y),p\>-\frac{1}{2}\<\s_1\s_1^*(x,y)p, p\>, \quad \bar{\Psi}(x,p):=\int_{\mR^m}\Psi(x,y,p)\mu^x(\dif y),
$$
and it holds that $(\sL^x_1\kappa)(x,y,p)=-(\Psi(x,y,p)-\bar{\Psi}(x,p))$.

\bp\label{poisequasoluchar}
Suppose that $(\mathbf{H}^1_{b_{1}, \s_{1}})$, $(\mathbf{H}^2_{b_{1}, \s_{1}})$, $(\mathbf{H}^1_{b_{2}, \s_{2}})$, $(\mathbf{H}^2_{b_{2}, \s_{2}})$ and $(\mathbf{H}^4_{b_{2}, \s_{2}})$ hold. Set
\be
\kappa(x,y,p)&:=&\int_0^\infty P_t^x[\Psi(x,\cdot,p)-\bar{\Psi}(x,p)](y)\dif t\no\\
&=&\int_0^\infty \left(\mE\Psi(x,Y_t^{x,y},p)-\bar{\Psi}(x,p)\right)\dif t.
\label{poisequasolu}
\ee
Then for any $x\in\overline{\cD(A)}$ and $p\in\mR^n$, $\kappa(x,\cdot,p)$ belongs to $C^2(\mR^m)$ and $\kappa(x,y,p)$ is the unique solution for Eq.(\ref{poisequa}). Moreover, it holds that
\ce
|\kappa(x,y,p)|\leq C(|p|+|p|^2)(1+|x|+|y|), \quad |\p_y \kappa(x,y,p)|\leq C(|p|+|p|^2).
\de
\ep
\begin{proof}
First of all, from $(\mathbf{H}^1_{b_{1}, \s_{1}})$ and $(\mathbf{H}^2_{b_1,\s_{1}})$, it follows that for $y_1, y_2\in\mR^m$
$$
|\Psi(x,y_1,p)-\Psi(x,y_2,p)|^2\leq (2L_{b_1,\s_1}|p|^2+C|p|^4)|y_1-y_2|^2,
$$
which together with the definition of $\mu^{x}$ and Lemma \ref{frozesti} implies that
\be
&&\left|P_{t}^{x}\left[\Psi(x,\cdot,p)-\bar{\Psi}(x,p)\right](y)\right|^2\no\\
&=&|\mE \Psi(x,Y_{t}^{x,y},p)-\bar{ \Psi}(x,p)|^2\no\\
&=&|\mE \Psi(x,Y_{t}^{x,y},p)-\int_{\mR^m}\Psi(x,z,p)\mu^{x}(\dif z)|^2\no\\
&=&|\mE \Psi(x,Y_{t}^{x,y},p)-\int_{\mR^m}\mE \Psi(x,Y_{t}^{x,z},p)\mu^{x}(\dif z)|^2\no\\
&\leq&\int_{\mR^m}\mE|\Psi(x,Y_{t}^{x,y},p)-\Psi(x,Y_{t}^{x,z},p)|^2\mu^{x}(\dif z)\no\\
&\leq&(2L_{b_1,\s_1}|p|^2+C|p|^4)\int_{\mR^m}|y-z|^2e^{-\a t}\mu^{x}(\dif z)\no\\
&\leq&C(2L_{b_1,\s_1}|p|^2+C|p|^4)e^{-\a t}(1+|x|+|y|)^2.
\label{semiesti}
\ee
Thus, it holds that
\be 
&&\int_{0}^{\infty}\left|P_{t}^{x}\left[\Psi(x,\cdot,p)-\bar{\Psi}(x,p)\right](y)\right|\dif t\no\\
&=&\int_{0}^{\infty}|\mE \Psi(x,Y_{t}^{x,y},p)-\bar{ \Psi}(x,p)|\dif t\no\\
&\leq&\int_{0}^{\infty}C(2^{1/2}L^{1/2}_{b_1,\s_1}|p|+C|p|^2)e^{-\frac{\a t}{2}}(1+|x|+|y|)\dif t\no\\
&=&\frac{2C}{\a}(2^{1/2}L^{1/2}_{b_1,\s_1}|p|+C|p|^2)(1+|x|+|y|).
\label{xfbo}
\ee
So the right side of (\ref{poisequasolu}) is well-defined. 

Next, we show that for any $x\in\overline{\cD(A)}$ and $p\in\mR^n$, $\kappa(x,\cdot,p)$ belongs to $C^2(\mR^m)$ and $\kappa(x,y,p)$ is the unique solution for Eq.(\ref{poisequa}). 

First of all, we study the regularity of $\Psi$. By $(\mathbf{H}^2_{b_{1}, \s_{1}})$, it holds that for any $x\in\overline{\cD(A)}$ and $p\in\mR^n$, $\Psi(x,\cdot,p)$ belongs to $C^2(\mR^m)$. Besides, note that $Y_{\cdot}^{x,y}$ satisfies Eq.(\ref{frozequa}) with $Y_0^{x,y}=y$. So, $(\mathbf{H}^{4}_{b_{2}, \s_{2}})$ assures the existence of $\partial_{y} Y_{t}^{x,y}, \partial_{yy} Y_{t}^{x,y}$.

Combining the above deduction, we have that for any $x\in\overline{\cD(A)}$ and $p\in\mR^n$, $\kappa(x,\cdot,p)$ belongs to $C^2(\mR^m)$. Then by acting the generator $\sL^x_1$ on $\kappa(x,y,p)$, it holds that 
\ce
(\sL^x_1\kappa)(x,y,p)&=&\int_{0}^{\infty}(\sL^x_1P_{t}^{x})\left[\Psi(x,\cdot,p)-\bar{\Psi}(x,p)\right](y)\dif t\no\\
&=&\int_{0}^{\infty}\frac{\dif P_{t}^{x}\left[\Psi(x,\cdot,p)-\bar{\Psi}(x,p)\right](y)}{\dif t}\dif t\no\\
&=& \lim_{t\rightarrow\infty}P_{t}^{x}\left[\Psi(x,\cdot,p)-\bar{\Psi}(x,p)\right](y)-\left[\Psi(x,\cdot,p)-\bar{\Psi}(x,p)\right](y)\no\\
&=&-\left[\Psi(x,\cdot,p)-\bar{\Psi}(x,p)\right](y),
\de
which yields that $\kappa(x,y,p)$ is a solution for Eq.(\ref{poisequa}). Moreover, based on $(\mathbf{H}^{4}_{b_{2}, \s_{2}})$, we know that the solutions of Eq.(\ref{poisequa}) are unique up to an additive constant. Thus, $\kappa(x,y,p)$ is the unique solution for Eq.(\ref{poisequa}). 

Finally, we establish the required estimates. By (\ref{xfbo}), we conclude that 
$$
|\kappa(x,y,p)|\leq C(|p|+|p|^2)(1+|x|+|y|).
$$
And from Lemma \ref{frozesti}, it follows that $|\p_{y}\kappa(x,y,p)|\leq C(|p|+|p|^2)$. The proof is complete.
\end{proof}

Next, we study  the relationship among $H_0(x, p), H_1(x, p)$ and $\bar{H}(x,p)$.

\bp\label{relah0h1barh}
Assume that $(\mathbf{H}^1_{b_{1}, \s_{1}})$, $(\mathbf{H}^2_{b_{1}, \s_{1}})$, $(\mathbf{H}^1_{b_{2}, \s_{2}})$-$(\mathbf{H}^4_{b_{2}, \s_{2}})$ hold. Then it holds that
\ce
H_0(x, p)\leq \bar{H}(x,p)\leq H_1(x, p), \quad x\in\overline{\cD(A)},\quad  p\in\mR^n.
\de
\ep
\begin{proof}
In order to prove that $H_0(x, p)\leq \bar{H}(x,p)$, by the definitions of $H_0(x, p)$ and $\bar{H}(x,p)$, we only need to choose a special $\xi$. 

By Lemma \ref{poisequasoluchar}, we know that $\kappa(x,y,p)$ belongs to $C^2(\mR^m)$ in $y$. But $\kappa(x,y,p)$ does not have a compact support in $y$. Taking $\rho\in C^\infty(\mR^m)$ such that $\rho(z)=1$ when $|z|\leq 1$ and $\rho(z)=0$ when $|z|\geq 2$. Then $\rho(z), \p_z\rho(z), \p_{zz}\rho(z)$ are uniformly bounded. For $l\in\mN$, set $\xi_l(y):=\rho(\frac y l)\kappa(x,y,p)$, and it holds that $\xi_l\in C_c^2(\mR^m)$ and
\ce
\sL^x_1\xi_l(y)&=&\frac 1 l\<b_2(x,y),\p_z\rho(\frac y l)\>\kappa(x,y,p)+\frac{1}{2l^2}{\rm tr}(\s_2\s_2^*(x,y)\p_{zz}\rho(\frac y l))\kappa(x,y,p)\\
&&+\frac 1 l {\rm tr}(\s_2\s_2^*(x,y)\p_z\rho(\frac y l)\p_y\kappa(x,y,p))+\rho(\frac y l)\sL^x_1\kappa(x,y,p).
\de

Next, in terms of the definition of $H_0(x, p)$ and Lemma \ref{poisequasoluchar}, we infer that
\ce
H_0(x, p)&=&\inf\limits _{\l \in \Lambda} H_0(x, p; \l)\\
&=&\inf\limits _{0<\t<1}\inf\limits _{\xi\in C_c^2(\mR^m)}\sup\limits_{y\in\mR^m}\bigg\{\<b_1(x,y),p\>-\frac{1}{2}\<\s_1\s^*_1(x,y)p,p\>+\sL^x_1\xi(y)\\
&&\qquad\quad-\t L_1\zeta(y)+\t L_2 I_B(y)\bigg\}\\
&\leq&\inf\limits _{0<\t<1}\limsup\limits_{l\rightarrow\infty}\sup\limits_{y\in\mR^m}\left\{\<b_1(x,y),p\>-\frac{1}{2}\<\s_1\s^*_1(x,y)p,p\>+\sL^x_1\xi_l(y)+\t L_2 I_B(y)\right\}\\
&\leq&\inf\limits _{0<\t<1}\sup\limits_{y\in\mR^m}\left\{\<b_1(x,y),p\>-\frac{1}{2}\<\s_1\s^*_1(x,y)p,p\>+\sL^x_1\kappa(x,y,p)+\t L_2 I_B(y)\right\}\\
&=&\<\bar{b}_1(x),p\>-\frac{1}{2}\<\bar{a}_1(x)p, p\>\\
&=&\bar{H}(x,p).
\de

Finally, by the similar deduction to that for $H_0(x, p)\leq \bar{H}(x,p)$, we obtain that $\bar{H}(x,p)\leq H_1(x, p)$. The proof is complete.
\end{proof}

\bp\label{compprin}
Assume that $(\mathbf{H}_{A})$, $(\mathbf{H}^1_{b_{1}, \s_{1}})$, $(\mathbf{H}^1_{b_{2}, \s_{2}})$-$(\mathbf{H}^2_{b_{2}, \s_{2}})$ hold and $h\in {\rm Lip}_b(\overline{\cD(A)})$. Suppose that $u, v$ are a viscosity subsolution and a viscosity supersolution for Eq.(\ref{hjblimi}), respectively. Then it holds that $u\leq v$.
\ep
\begin{proof}
First, by Lemma \ref{barb1bara1lip}, it holds that $\bar{H}(x,p)$ is continuous in $(x,p)$. 

In the following, it holds that for any $R>0$ and $|p_1|, |p_2|\leq R$
\ce
&&|\bar{H}(x,p_1)-\bar{H}(x,p_2)|\\
&=&\left|\<\bar{b}_1(x),p_1\>-\frac{1}{2}\<\bar{a}_1(x)p_1, p_1\>-\<\bar{b}_1(x),p_2\>+\frac{1}{2}\<\bar{a}_1(x)p_2, p_2\>\right|\\
&\leq&L^{1/2}_{b_{1},\s_{1}}|p_1-p_2|+\frac{1}{2}L_{b_{1},\s_{1}}(|p_1|+|p_2|)|p_1-p_2|\\
&\leq&(L^{1/2}_{b_{1},\s_{1}}+L_{b_{1},\s_{1}}R)|p_1-p_2|.
\de
Define $\varpi_R: [0,\infty)\to[0,\infty)$ by $\varpi_R(r)=(L^{1/2}_{b_{1},\s_{1}}+L_{b_{1},\s_{1}}R)r$, and then $\varpi_R$  is continuous increasing and $\varpi_R(0)=0$. Moreover,
\ce
|\bar{H}(x,p_1)-\bar{H}(x,p_2)|\leq \varpi_R(|p_1-p_2|).
\de

Next, we verify that $u(t,x)$ or $v(t,x)$ is Lipschitz continuous in $x$. Let us observe Eq.(\ref{ratefuncequa}). For any $x_1, x_2\in\overline{\cD(A)}$ and $t\in[0,T]$, it holds that
\ce
|X_{x_1}^z(t)-X_{x_2}^z(t)|^2&=&|x_1-x_2|^2-2\int_0^t\<X_{x_1}^z(r)-X_{x_2}^z(r), \dif (K_{x_1}^z(r)-K_{x_2}^z(r))\>\\
&&+2\int_0^t\<X_{x_1}^z(r)-X_{x_2}^z(r),\bar{b}_1(X_{x_1}^z(r))-\bar{b}_1(X_{x_2}^z(r))\>\dif r\\
&&+2\int_0^t\<X_{x_1}^z(r)-X_{x_2}^z(r),(\bar{\s}_1(X_{x_1}^z(r))-\bar{\s}_1(X_{x_2}^z(r)))z(r)\>\dif r.
\de
By Lemma \ref{equi}, we know that
$$
-2\int_0^t\<X_{x_1}^z(r)-X_{x_2}^z(r), \dif (K_{x_1}^z(r)-K_{x_2}^z(r))\>\leq 0.
$$
Lemma \ref{barb1bara1lip} implies that
\ce
2\int_0^t\<X_{x_1}^z(r)-X_{x_2}^z(r),\bar{b}_1(X_{x_1}^z(r))-\bar{b}_1(X_{x_2}^z(r))\>\dif r\leq 2C\int_0^t|X_{x_1}^z(r)-X_{x_2}^z(r)|^2\dif r,
\de
and
\ce
&&2\int_0^t\<X_{x_1}^z(r)-X_{x_2}^z(r),(\bar{\s}_1(X_{x_1}^z(r))-\bar{\s}_1(X_{x_2}^z(r)))z(r)\>\dif r\\
&\leq&2\sup\limits_{r\in[0,t]}|X_{x_1}^z(r)-X_{x_2}^z(r)|\int_0^t|\bar{\s}_1(X_{x_1}^z(r))-\bar{\s}_1(X_{x_2}^z(r))||z(r)|\dif r\\
&\leq&\frac{1}{2}\sup\limits_{r\in[0,t]}|X_{x_1}^z(r)-X_{x_2}^z(r)|^2+2\left(\int_0^t|\bar{\s}_1(X_{x_1}^z(r))-\bar{\s}_1(X_{x_2}^z(r))||z(r)|\dif r\right)^2\\
&\leq&\frac{1}{2}\sup\limits_{r\in[0,t]}|X_{x_1}^z(r)-X_{x_2}^z(r)|^2+2\int_0^t|\bar{\s}_1(X_{x_1}^z(r))-\bar{\s}_1(X_{x_2}^z(r))|^2\dif r\int_0^t|z(r)|^2\dif r\\
&\leq&\frac{1}{2}\sup\limits_{r\in[0,t]}|X_{x_1}^z(r)-X_{x_2}^z(r)|^2+2C^2\left(\int_0^T|z(r)|^2\dif r\right)\int_0^t|X_{x_1}^z(r)-X_{x_2}^z(r)|^2\dif r.
\de
These deductions yield that
\ce
&&\sup\limits_{r\in[0,t]}|X_{x_1}^z(r)-X_{x_2}^z(r)|^2\\
&\leq&2|x_1-x_2|^2+2\left(C+\left(\int_0^T|z(r)|^2\dif r\right)\right)\int_0^t\sup\limits_{r\in[0,\nu]}|X_{x_1}^z(r)-X_{x_2}^z(r)|^2\dif \nu.
\de
So, the Gronwall inequality implies that
\ce
\sup\limits_{r\in[0,t]}|X_{x_1}^z(r)-X_{x_2}^z(r)|^2\leq 2|x_1-x_2|^2 e^{CT}.
\de
In terms of this and $h\in {\rm Lip}_b(\overline{\cD(A)})$, we infer that $u^h_0(t,x_0)$ defined by (\ref{hjblimisolu}) is Lipschitz continuous in $x_0$. Besides, Theorem 4 in \cite{dT} means that $u^h_0(t,x_0)$ is a viscosity solution for Eq.(\ref{hjblimi}). From this, it follows that $u(t,x)$ or $v(t,x)$ is Lipschitz continuous in $x$.

Finally, collecting the above deduction, by Theorem \ref{compprinth} in Appendix we conclude that $u\leq v$.
\end{proof}

Now, combining Proposition \ref{relah0h1barh} with Proposition \ref{compprin}, by Definition \ref{viscsoludefi} we conclude that a comparison principle between the viscosity subsolutions to Eq.(\ref{conchohhjb}) and the viscosity supersolutions to Eq.(\ref{conch1hhjb}) holds. Then Proposition \ref{concsupersubsolu} and Lemma \ref{uehlimi} yield the following result.

\bp\label{ueghuohhlip}
Suppose that $(\mathbf{H}_{A})$, $(\mathbf{H}^1_{b_{1}, \s_{1}})$, $(\mathbf{H}^2_{b_{1}, \s_{1}})$, $(\mathbf{H}^1_{b_{2}, \s_{2}})$-$(\mathbf{H}^4_{b_{2}, \s_{2}})$ hold and $h\in {\rm Lip}_b(\overline{\cD(A)})$. If 
\ce
\lim\limits_{\e\rightarrow0}\frac{\g}{\e}=0,
\de
then 
$$
\lim\limits_{\e\rightarrow 0}\sup\limits_{t\in[0,T]}\sup\limits_{(x,y)\in K \times \tilde{K}}|u_{\e,\g}^h(t,x,y)-u_0^h(t,x)|=0, \quad \forall K \times \tilde{K}\in\mathcal{Q},
$$
where $u_0^h$ is the unique viscosity solution of Eq.(\ref{hjblimi}).
\ep
\begin{proof}
Firstly, Proposition \ref{concsupersubsolu}, Proposition \ref{relah0h1barh}, Proposition \ref{compprin} and Lemma \ref{uehlimi} imply that $u^h:=\bar{u}^h=\underline{u}^h$ and 
$$
\lim\limits_{\e\rightarrow 0}\sup\limits_{t\in[0,T]}\sup\limits_{(x,y)\in K \times \tilde{K}}|u_{\e,\g}^h(t,x,y)-u^h(t,x)|=0, \quad \forall K \times \tilde{K}\in\mathcal{Q}.
$$

Next, we prove that $u^h$ is the unique viscosity solution of Eq.(\ref{hjblimi}). On the one hand, Proposition \ref{concsupersubsolu} and Proposition \ref{relah0h1barh} infer that $\bar{u}^h$ and $\underline{u}^h$ are a viscosity subsolution and a viscosity supersolution for Eq.(\ref{hjblimi}), respectively. Note that $u^h=\bar{u}^h=\underline{u}^h$. Thus, $u^h$ is a viscosity solution for Eq.(\ref{hjblimi}). On the other hand, Proposition \ref{compprin} assures the uniqueness of viscosity solutions for Eq.(\ref{hjblimi}). So, we conclude that $u^h$ is the unique viscosity solution of Eq.(\ref{hjblimi}). The proof is complete.
\end{proof}

In order to obtain that $u_{\e,\g}^h$ converges to $u_0^h$ for $h\in C_b(\overline{\cD(A)})$, we need the following exponential estimate.

\bp\label{expoesti}
Suppose that $(\mathbf{H}_{A})$, $(\mathbf{H}^1_{b_{1}, \s_{1}})$ and $(\mathbf{H}^1_{b_{2}, \s_{2}})$ hold. Then there exist two constants $C_1, C_2>0$ independent of $\e, \g$ such that
\ce
\mE\sup\limits_{t\in[0,T]}\exp\left\{\frac{C_1}{\e}\(\log(e+|X^{\e,\g}_t|^2)\)^2\right\}\leq C_2 e^{\frac{2(\log(e+|x_0|^2))^2}{\e}}.
\de
\ep
\begin{proof}
Let $\vartheta(x):=\(\log(e+|x|^2)\)^2$ for $x\in\mR^n$ and $\tau_R:=\inf\{t\leq T: |X^{\e,\g}_t|>R\}$. Taking any $\nu>0, \iota>0$ and applying the It\^o formula to $e^{\frac{\nu}{\e} e^{-\iota t}\vartheta(X^{\e,\g}_t)}$, we obtain that
\ce
&&e^{\frac{\nu}{\e} e^{-\iota (t\land \tau_R)}\vartheta(X^{\e,\g}_{t\land \tau_R})}\no\\
&=&e^{\frac{\nu}{\e}\vartheta(x_0)}+\int_0^{t\land \tau_R}e^{\frac{\nu}{\e} e^{-\iota s}\vartheta(X^{\e,\g}_s)}\frac{\nu}{\e} e^{-\iota s}(-\iota)\vartheta(X^{\e,\g}_s)\dif s\no\\
&&-4\int_0^{t\land \tau_R}e^{\frac{\nu}{\e} e^{-\iota s}\vartheta(X^{\e,\g}_s)}\frac{\nu}{\e} e^{-\iota s}\frac{\log(e+|X^{\e,\g}_s|^2)}{e+|X^{\e,\g}_s|^2}\<X^{\e,\g}_s, \dif K^{1,\e,\g}_s\>\no\\
&&+4\int_0^{t\land \tau_R}e^{\frac{\nu}{\e} e^{-\iota s}\vartheta(X^{\e,\g}_s)}\frac{\nu}{\e} e^{-\iota s}\frac{\log(e+|X^{\e,\g}_s|^2)}{e+|X^{\e,\g}_s|^2}\<X^{\e,\g}_s, b_{1}(X_{s}^{\e,\g},Y_{s}^{\e,\g})\>\dif s\no\\
&&+4\int_0^{t\land \tau_R}e^{\frac{\nu}{\e} e^{-\iota s}\vartheta(X^{\e,\g}_s)}\frac{\nu}{\e} e^{-\iota s}\frac{\log(e+|X^{\e,\g}_s|^2)}{e+|X^{\e,\g}_s|^2}\sqrt \e\<X^{\e,\g}_s, \s_{1}(X_{s}^{\e,\g},Y_{s}^{\e,\g})\dif W^1_{s}\>\no\\
&&+8\int_0^{t\land \tau_R}e^{\frac{\nu}{\e} e^{-\iota s}\vartheta(X^{\e,\g}_s)}\frac{\nu^2}{\e} e^{-2\iota s}\frac{(\log(e+|X^{\e,\g}_s|^2))^2}{(e+|X^{\e,\g}_s|^2)^2}|\s^*_{1}(X_{s}^{\e,\g},Y_{s}^{\e,\g})X^{\e,\g}_s|^2\dif s \no\\
&&+4\int_0^{t\land \tau_R}e^{\frac{\nu}{\e} e^{-\iota s}\vartheta(X^{\e,\g}_s)}\frac{\nu}{\e} e^{-\iota s}\frac{\e}{(e+|X^{\e,\g}_s|^2)^2}|\s^*_{1}(X_{s}^{\e,\g},Y_{s}^{\e,\g})X^{\e,\g}_s|^2\dif s\no \\
&&+2\int_0^{t\land \tau_R}e^{\frac{\nu}{\e} e^{-\iota s}\vartheta(X^{\e,\g}_s)}\frac{\nu}{\e} e^{-\iota s}\e\log(e+|X^{\e,\g}_s|^2)\bigg[\frac{{\rm tr}(\s\s^*_{1}(X_{s}^{\e,\g},Y_{s}^{\e,\g}))}{e+|X^{\e,\g}_s|^2}\\
&&\qquad -\frac{2|\s^*_{1}(X_{s}^{\e,\g},Y_{s}^{\e,\g})X^{\e,\g}_s|^2}{(e+|X^{\e,\g}_s|^2)^2}\bigg].
\de

Next, by Lemma \ref{equi}, it holds that for any $\varUpsilon\in A(0)$,
\ce
\<X^{\e,\g}_s-0, \dif K^{1,\e,\g}_s-\varUpsilon\dif s\>\geq 0,
\de
and
\ce
-\<X^{\e,\g}_s,\dif K^{1,\e,\g}_s\>\leq -\<X^{\e,\g}_s,\varUpsilon\>\dif s\leq |X^{\e,\g}_s||\varUpsilon|\dif s\leq |\varUpsilon|(e+|X^{\e,\g}_s|^2)\dif s.
\de
Moreover, $(\mathbf{H}^1_{b_{1}, \s_{1}})$ implies that
\ce
&&\<X^{\e,\g}_s, b_{1}(X_{s}^{\e,\g},Y_{s}^{\e,\g})\>\leq |X^{\e,\g}_s||b_{1}(X_{s}^{\e,\g},Y_{s}^{\e,\g})|\leq L^{1/2}_{b_1,\s_1}(e+|X^{\e,\g}_s|^2),\\
&&|\s^*_{1}(X_{s}^{\e,\g},Y_{s}^{\e,\g})X^{\e,\g}_s|^2\leq \|\s^*_{1}(X_{s}^{\e,\g},Y_{s}^{\e,\g})\|^2|X^{\e,\g}_s|^2\leq L_{b_1,\s_1}(e+|X^{\e,\g}_s|^2).
\de
Collecting the above deduction, we have that
\be
&&e^{\frac{\nu}{\e} e^{-\iota (t\land \tau_R)}\vartheta(X^{\e,\g}_{t\land \tau_R})}\no\\
&\leq& e^{\frac{\nu}{\e}\vartheta(x_0)}-\iota\int_0^{t\land \tau_R}e^{\frac{\nu}{\e} e^{-\iota s}\vartheta(X^{\e,\g}_s)}\frac{\nu}{\e} e^{-\iota s}\vartheta(X^{\e,\g}_s)\dif s\no\\
&&+4\int_0^{t\land \tau_R}e^{\frac{\nu}{\e} e^{-\iota s}\vartheta(X^{\e,\g}_s)}\frac{\nu}{\e} e^{-\iota s}\frac{\log(e+|X^{\e,\g}_s|^2)}{e+|X^{\e,\g}_s|^2}\sqrt \e\<X^{\e,\g}_s, \s_{1}(X_{s}^{\e,\g},Y_{s}^{\e,\g})\dif W^1_{s}\>\no\\
&&+(4|\varUpsilon|+4L^{1/2}_{b_1,\s_1}+8\nu L_{b_1,\s_1}+10 L_{b_1,\s_1})\int_0^{t\land \tau_R}e^{\frac{\nu}{\e} e^{-\iota s}\vartheta(X^{\e,\g}_s)}\frac{\nu}{\e} e^{-\iota s}\vartheta(X^{\e,\g}_s)\dif s.
\label{expoito}
\ee
By taking the expectation on two sides of (\ref{expoito}), it follows that
\ce
&&\mE e^{\frac{\nu}{\e} e^{-\iota (t\land \tau_R)}\vartheta(X^{\e,\g}_{t\land \tau_R})}\\
&\leq& e^{\frac{\nu}{\e}\vartheta(x_0)}-\iota\mE\int_0^{t\land \tau_R}e^{\frac{\nu}{\e} e^{-\iota s}\vartheta(X^{\e,\g}_s)}\frac{\nu}{\e} e^{-\iota s}\vartheta(X^{\e,\g}_s)\dif s\\
&&+(4|\varUpsilon|+4L^{1/2}_{b_1,\s_1}+8\nu L_{b_1,\s_1}+10 L_{b_1,\s_1})\mE\int_0^{t\land \tau_R}e^{\frac{\nu}{\e} e^{-\iota s}\vartheta(X^{\e,\g}_s)}\frac{\nu}{\e} e^{-\iota s}\vartheta(X^{\e,\g}_s)\dif s.
\de
Letting $\iota=4|\varUpsilon|+4L^{1/2}_{b_1,\s_1}+8\nu L_{b_1,\s_1}+10 L_{b_1,\s_1}+1$, one can infer that
\ce
\mE e^{\frac{\nu}{\e} e^{-\iota (t\land \tau_R)}\vartheta(X^{\e,\g}_{t\land \tau_R})}+\mE\int_0^{t\land \tau_R}e^{\frac{\nu}{\e} e^{-\iota s}\vartheta(X^{\e,\g}_s)}\frac{\nu}{\e} e^{-\iota s}\vartheta(X^{\e,\g}_s)\dif s\leq e^{\frac{\nu}{\e}\vartheta(x_0)}.
\de
We take $\nu=2$ and $\iota=4|\varUpsilon|+4L^{1/2}_{b_1,\s_1}+16 L_{b_1,\s_1}+10 L_{b_1,\s_1}+1$ and get that
\ce
\mE e^{\frac{2}{\e} e^{-\iota (t\land \tau_R)}\vartheta(X^{\e,\g}_{t\land \tau_R})}+\mE\int_0^{t\land \tau_R}e^{\frac{2}{\e} e^{-\iota s}\vartheta(X^{\e,\g}_s)}\frac{2}{\e} e^{-\iota s}\vartheta(X^{\e,\g}_s)\dif s\leq e^{\frac{2}{\e}\vartheta(x_0)},
\de
and furthermore
\ce
\mE\int_0^{t\land \tau_R}e^{\frac{1}{\e} e^{-\iota s}\vartheta(X^{\e,\g}_s)}\frac{1}{\e} e^{-\iota s}\vartheta(X^{\e,\g}_s)\dif s\leq \mE\int_0^{t\land \tau_R}e^{\frac{2}{\e} e^{-\iota s}\vartheta(X^{\e,\g}_s)}\frac{2}{\e} e^{-\iota s}\vartheta(X^{\e,\g}_s)\dif s\leq e^{\frac{2}{\e}\vartheta(x_0)}.
\de

In the following, based on (\ref{expoito}) with $\nu=1$ and $\iota=4|\varUpsilon|+4L^{1/2}_{b_1,\s_1}+8L_{b_1,\s_1}+10 L_{b_1,\s_1}+1$, the BDG inequality and the H\"older inequality, it holds that
\ce
&&\mE\sup\limits_{t\in[0,T]}e^{\frac{1}{\e} e^{-\iota (t\land \tau_R)}\vartheta(X^{\e,\g}_{t\land \tau_R})}\\
&\leq&e^{\frac{1}{\e}\vartheta(x_0)}+24L^{1/2}_{b_1,\s_1}\left(\mE\int_0^{T\land \tau_R}e^{\frac{2}{\e} e^{-\iota s}\vartheta(X^{\e,\g}_s)}\frac{1}{\e} e^{-2\iota s}\vartheta(X^{\e,\g}_s)\dif s\right)^{1/2}\\
&&+(4|\varUpsilon|+4L^{1/2}_{b_1,\s_1}+8L_{b_1,\s_1}+10 L_{b_1,\s_1})\mE\int_0^{T\land \tau_R}e^{\frac{1}{\e} e^{-\iota s}\vartheta(X^{\e,\g}_s)}\frac{1}{\e} e^{-\iota s}\vartheta(X^{\e,\g}_s)\dif s\\
&\leq&e^{\frac{1}{\e}\vartheta(x_0)}+24L^{1/2}_{b_1,\s_1}e^{\frac{1}{\e}\vartheta(x_0)}+(4|\varUpsilon|+4L^{1/2}_{b_1,\s_1}+18 L_{b_1,\s_1})e^{\frac{2}{\e}\vartheta(x_0)}.
\de
Letting $R\rightarrow\infty$, we conclude that
\ce
\mE\sup\limits_{t\in[0,T]}e^{\frac{1}{\e} e^{-\iota t}\vartheta(X^{\e,\g}_{t})}\leq (1+4|\varUpsilon|+28L^{1/2}_{b_1,\s_1}+18 L_{b_1,\s_1})e^{\frac{2}{\e}\vartheta(x_0)}.
\de

Finally, taking $C_1=e^{-\iota T}$ and $C_2=1+4|\varUpsilon|+28L^{1/2}_{b_1,\s_1}+18 L_{b_1,\s_1}$, we obtain the required estimate. The proof is complete.
\end{proof}

We also need the following convergence result.

\bp\label{ueghkh}
Suppose that $(\mathbf{H}_{A})$, $(\mathbf{H}^1_{b_{1}, \s_{1}})$ and $(\mathbf{H}^1_{b_{2}, \s_{2}})$ hold. Assume that $h_k, h\in C_b(\overline{\cD(A)})$ for $k\in\mN$ and $h_k$ converges to $h$ in $C_b(\overline{\cD(A)})$. Then for any $R>0$, any compact set $\tilde{K}\subset\mR^m$ and $\eta>0$ there exist a $k_0\in\mN$ and a $\e_0\in (0,1)$ such that if $k\geq k_0$ and $\e<\e_0$
\be\label{ueghkhes}
\sup\limits_{x\in B(0,R)\cap\overline{\cD(A)}, y\in \tilde{K}}|u^{h_k}_{\e,\g}(t,x,y)-u^{h}_{\e,\g}(t,x,y)|\leq \eta,
\ee
where $B(0,R):=\{x\in\mR^n: |x|\leq R\}$.
\ep
\begin{proof}
Let $x\in B(0,R)\cap\overline{\cD(A)}, y\in \tilde{K}$ and $R_1>R$. We take $\Omega_1:=\{\omega: |X^{\e,\g}_t(\omega)|>R_1\}$. So, by Proposition \ref{expoesti}, we choose $R_1$ such that
\be
\mP(\Omega_1)\leq e^{-\frac{3M}{\e}},
\label{omeges}
\ee
where $M\geq\sup\limits_{k\in \{0\}\cup\mN}\sup\limits_{x\in\overline{\cD(A)}}|h_k(x)|$ is a constant and $h_0:=h$.

Besides, we note that $h_k$ converges to $h$ in $C_b(\overline{\cD(A)})$. Thus, for any $\eta>0$ there exists a $k_0\in\mN$ such that for $k\geq k_0$
\be
\sup\limits_{x\in B(0,R_1)\cap\overline{\cD(A)}}|h_k(x)-h(x)|\leq \eta.
\label{hkhes}
\ee

In the following, we know that
\ce
&&\mE\exp\left(-\frac{h_k(X^{\e,\g}_t)}{\e}\right)\\
&=&\mE\exp\left(-\frac{h(X^{\e,\g}_t)}{\e}\right)+\mE\left[\exp\left(-\frac{h_k(X^{\e,\g}_t)}{\e}\right)-\exp\left(-\frac{h(X^{\e,\g}_t)}{\e}\right)\right]I_{\Omega\setminus\Omega_1}\\
&&+\mE\left[\exp\left(-\frac{h_k(X^{\e,\g}_t)}{\e}\right)-\exp\left(-\frac{h(X^{\e,\g}_t)}{\e}\right)\right]I_{\Omega_1}.
\de
By (\ref{omeges}) and (\ref{hkhes}), it holds that
\ce
&&\mE\exp\left(-\frac{h(X^{\e,\g}_t)}{\e}\right)+\mE\exp\left(-\frac{h(X^{\e,\g}_t)}{\e}\right)(e^{-\frac{\eta}{\e}}-1)-2e^{-\frac{2M}{\e}}\\
&\leq&\mE\exp\left(-\frac{h_k(X^{\e,\g}_t)}{\e}\right)\\
&\leq&\mE\exp\left(-\frac{h(X^{\e,\g}_t)}{\e}\right)+\mE\exp\left(-\frac{h(X^{\e,\g}_t)}{\e}\right)(e^{\frac{\eta}{\e}}-1)+2e^{-\frac{2M}{\e}},
\de
and
\ce
&&\e\log\mE\exp\left(-\frac{h(X^{\e,\g}_t)}{\e}\right)-\eta+\e\log\left[1-\frac{2e^{-\frac{2M}{\e}}}{\mE\exp\left(-\frac{h(X^{\e,\g}_t)}{\e}\right) e^{-\frac{\eta}{\e}}}\right]\\
&\leq&\e\log\mE\exp\left(-\frac{h_k(X^{\e,\g}_t)}{\e}\right)\\
&\leq&\e\log\mE\exp\left(-\frac{h(X^{\e,\g}_t)}{\e}\right)+\eta+\e\log\left[1+\frac{2e^{-\frac{2M}{\e}}}{\mE\exp\left(-\frac{h(X^{\e,\g}_t)}{\e}\right) e^{-\frac{\eta}{\e}}}\right].
\de
Note that for small $s>0$
$$
-2s\leq \log(1-s), \quad \log(1+s)\leq s.
$$
Thus, we infer that
\ce
&&\e\log\mE\exp\left(-\frac{h(X^{\e,\g}_t)}{\e}\right)-\eta-4\e e^{-\frac{M-\eta}{\e}}\\
&\leq&\e\log\mE\exp\left(-\frac{h_k(X^{\e,\g}_t)}{\e}\right)\\
&\leq&\e\log\mE\exp\left(-\frac{h(X^{\e,\g}_t)}{\e}\right)+\eta+2\e e^{-\frac{M-\eta}{\e}}.
\de
By taking $\eta<M$ and enough small $\e$, (\ref{ueghkhes}) holds, which completes the proof.
\end{proof}

Now, we give the main result in this subsection.

\bt\label{ueghuoh}
Suppose that $(\mathbf{H}_{A})$, $(\mathbf{H}^1_{b_{1}, \s_{1}})$, $(\mathbf{H}^2_{b_{1}, \s_{1}})$, $(\mathbf{H}^1_{b_{2}, \s_{2}})$-$(\mathbf{H}^4_{b_{2}, \s_{2}})$ hold and $h\in C_b(\overline{\cD(A)})$. If 
\ce
\lim\limits_{\e\rightarrow0}\frac{\g}{\e}=0,
\de
then 
$$
\lim\limits_{\e\rightarrow 0}\sup\limits_{t\in[0,T]}\sup\limits_{(x,y)\in K \times \tilde{K}}|u_{\e,\g}^h(t,x,y)-u_0^h(t,x)|=0, \quad \forall K \times \tilde{K}\in\mathcal{Q}.
$$
\et
\begin{proof}
First of all, for any $h\in C_b(\overline{\cD(A)})$, there exists a sequence $\{h_k\}\subset{\rm Lip}_b(\overline{\cD(A)})$ such that $h_k$ converges to $h$ in $C_b(\overline{\cD(A)})$. By Proposition \ref{ueghkh},  for any $R>0$, any compact set $\tilde{K}\subset\mR^m$ and $\eta>0$ there exist a $k_0\in\mN$ and a $\e_0\in (0,1)$ such that if $k_1, k_2\geq k_0$ and $\e<\e_0$
\ce
\sup\limits_{x\in B(0,R)\cap\overline{\cD(A)}, y\in \tilde{K}}|u^{h_{k_1}}_{\e,\g}(t,x,y)-u^{h_{k_2}}_{\e,\g}(t,x,y)|\leq 2\eta.
\de
By taking the limit on two sides of the above inequality as $\e$ tends to $0$, Proposition \ref{ueghuohhlip} implies that
\ce
\sup\limits_{x\in B(0,R)\cap\overline{\cD(A)}}|u^{h_{k_1}}_0(t,x)-u^{h_{k_2}}_{0}(t,x)|\leq 2\eta,
\de
which yields that $\lim\limits_{k\rightarrow\infty}u^{h_k}_0$ exists uniformly on compact subsets of $\overline{\cD(A)}$.

Besides, again by Proposition \ref{ueghkh} and Proposition \ref{ueghuohhlip}, it holds that for $\e<\e_0$
\ce
\sup\limits_{x\in B(0,R)\cap\overline{\cD(A)}, y\in \tilde{K}}|u_0^{h_{k_0}}(t,x)-u^{h}_{\e,\g}(t,x,y)|\leq 2\eta.
\de
As $k_0\rightarrow \infty$, by the above inequality we conclude that for $\e<\e_0$
\ce
\sup\limits_{x\in B(0,R)\cap\overline{\cD(A)}, y\in \tilde{K}}|\lim\limits_{k\rightarrow\infty}u^{h_k}_0(t,x)-u^{h}_{\e,\g}(t,x,y)|\leq 2\eta,
\de
which implies that
\ce
u_0^h(t,x)=\lim\limits_{k\rightarrow\infty}u^{h_k}_0(t,x)=\lim\limits_{\e\rightarrow 0}u^{h}_{\e,\g}(t,x,y),
\de
and
$$
\lim\limits_{\e\rightarrow 0}\sup\limits_{t\in[0,T]}\sup\limits_{(x,y)\in K \times \tilde{K}}|u_{\e,\g}^h(t,x,y)-u_0^h(t,x)|=0, \quad \forall K \times \tilde{K}\in\mathcal{Q}.
$$
The proof is complete.
\end{proof}

\subsection{Exponential tightness of $\{X^{\e,\g}_t, \e, \g>0\}$}

\bt\label{xegexpotigh}
Suppose that $(\mathbf{H}_{A})$, $(\mathbf{H}^1_{b_{1}, \s_{1}})$ and $(\mathbf{H}^1_{b_{2}, \s_{2}})$ hold. If for every $r>0$, there exists a compact subset $\Gamma_r\subset\overline{\cD(A)}$ such that for any $t\in[0,T]$
$$
\limsup\limits _{\e\rightarrow 0}\e \log \mP\left(X^{\e,\g}_t\notin \Gamma_r\right) \leq -r.
$$
\et
\begin{proof}
First of all, note that $\overline{\cD(A)}\subset\mR^n$. Thus, we take $\Gamma_r=\{x\in\overline{\cD(A)}: |x|\leq M\}$, where $M>0$ is determined later. By the Chebyshev inequality and Proposition \ref{expoesti}, it holds that for any $t\in[0,T]$
\ce
\mP\left(X^{\e,\g}_t\notin \Gamma_r\right)&=&\mP\left(|X^{\e,\g}_t|>M\right)\leq \mP\left(\sup\limits_{t\in[0,T]}|X^{\e,\g}_t|>M\right)\\
&=&\mP\left(\sup\limits_{t\in[0,T]}e^{\frac{C_1}{\e}\vartheta(X^{\e,\g}_{t})}>e^{\frac{C_1}{\e}\vartheta(M)}\right)\\
&\leq&\mE\sup\limits_{t\in[0,T]}e^{\frac{C_1}{\e}\vartheta(X^{\e,\g}_{t})}\times e^{-\frac{C_1}{\e}\vartheta(M)}\\
&\leq&C_2 e^{\frac{2(\log(e+|x_0|^2))^2}{\e}}\times e^{-\frac{C_1}{\e}\vartheta(M)}.
\de
In terms of the above inequality, we conclude that
\ce
\limsup\limits _{\e\rightarrow 0}\e \log \mP\left(X^{\e,\g}_t\notin \Gamma_r\right)\leq 2(\log(e+|x_0|^2))^2-C_1(\log(e+M^2))^2.
\de
By taking $M\geq \sqrt{\exp\{(r+2\vartheta(x_0))^{1/2} C_1^{-1/2}\}-e}$, it follows that
\ce
\limsup\limits _{\e\rightarrow 0}\e \log \mP\left(X^{\e,\g}_t\notin \Gamma_r\right)\leq -r,
\de
which completes the proof.
\end{proof}

Now, it is the position to prove Theorem \ref{ldptime}.

{\bf Proof of Theorem \ref{ldptime}.}  By Theorem \ref{ueghuoh} and Theorem \ref{xegexpotigh}, the Bryc formula (cf. \cite[Proposition 3.8]{fk}) yields that $\{X^{\e,\g}_t, \e, \g>0\}$ for any $t\in[0,T]$ satisfies the LDP with the speed $\e^{-1}$ and the good rate function $I$ given by
\ce
I(x; x_0,t)=\sup\limits_{h\in C_b(\overline{\cD(A)})}\left\{u_0^h(t, x_0)-h(x)\right\}.
\de
Based on Theorem 4 in \cite{dT} and the Varadhan lemma (cf. \cite[Proposition 3.8]{fk}), we conclude that
\ce
I(x; x_0, t)=\inf\left\{\frac{1}{2}\int_0^t|z(s)|^2\dif s: z\in L^2([0,T],\mR^n), s.t. X^z_{x_0}(t)=x\right\}.
\de
The proof is complete.

\section{An example}\label{exam}

In this section, we give an example to explain our result.

\bx
Consider the following slow-fast system on $\mR \times \mR$:
\be\left\{\begin{array}{l}
\dif X_{t}^{\e,\g}\in -\p\psi(X_{t}^{\e,\g})\dif t+\e\left[r-\frac{1}{2}\cos^2(Y_{t}^{\e,\g})\right]\dif t+\sqrt{\e}\cos(Y_{t}^{\e,\g})\dif W^1_{t},\\
X_{0}^{\e,\g}=x_0\in\overline{\cD(\p\psi)},\quad  0\leq t\leq T,\\
\dif Y_{t}^{\e,\g}=\frac{1}{\g}(\varsigma-\frac{1}{2}Y_{t}^{\e,\g})\dif t+\frac{\nu}{\sqrt{\g}}\dif W^2_{t},\\
Y_{0}^{\e,\g}=y_0,\quad  0\leq t\leq T,
\end{array}
\right.
\label{exameq}
\ee
where $\psi$ is a lower semicontinuous convex function, $\p\psi$ is a multivalued maximal monotone operator, $0\in{\rm Int}(\cD(\p\psi))$ and $r>0, \varsigma\in\mR, \nu>0$ are three constants. 

It is easy to see that
\ce
b_1(x,y)=r-\frac{1}{2}\cos^2(y), \quad \s_1(x,y)=\cos(y), \quad b_2(x,y)=\varsigma-\frac{1}{2}y, \quad \s_2(x,y)=\nu.
\de
By some computation, we know that $b_1, \s_1$ satisfy $(\mathbf{H}^{1}_{b_{1}, \s_{1}})$ and $(\mathbf{H}^{2}_{b_{1}, \s_{1}})$. 

Next, it holds that for $x_{i}\in\overline{\cD(\p\psi)}$, $y_{i}\in\mR$, $i=1, 2$
\ce
|b_2(x_1,y_1)-b_2(x_2,y_2)|^{2}+|\s_{2}(x_{1},y_{1})-\s_{2}(x_{2},y_{2})|^{2}=\frac{1}{4}|y_1-y_2|^2,
\de
and for $x\in\overline{\cD(\p\psi)}$, $y_{i}\in\mR$, $i=1, 2$
\ce
2(y_1-y_2)(b_2(x,y_1)-b_2(x,y_2))+|\s_2(x,y_1)-\s_2(x,y_2)|^2\leq -|y_1-y_2|^2,
\de
where $L_{b_2,\s_2}=\frac{1}{4}$ and $\b=1>2\times\frac{1}{4}$. So, $(\mathbf{H}^{1}_{b_{2}, \s_{2}})$ and $(\mathbf{H}^{2}_{b_{2}, \s_{2}})$ hold. And the following SDE 
\ce\left\{\begin{array}{l}
\dif Y_{t}=(\varsigma-\frac{1}{2}Y_{t})\dif t+\nu\dif W^2_{t},\\
Y_{0}=y_0,\quad  t\geq 0,
\end{array}
\right.
\de
has the unique invariant probability measure $\mu$.

In the following, set $\zeta(y):=|\frac{1}{2}y-\varsigma|^{\frac{3}{2}}$ and $\zeta$ has compact finite level sets. Moreover, 
\ce
(\sL^x_1\zeta)(y)&=&(\varsigma-\frac{1}{2}y)\p_y\zeta(y)+\frac{1}{2}\nu^2\p^2_{yy}\zeta(y)\\
&=&-\frac{3}{4}\zeta(y)+\frac{3\nu^2}{32}|\frac{1}{2}y-\varsigma|^{-\frac{1}{2}}.
\de
There exists a large enough $R>0$ such that 
\ce
(\sL^x_1\zeta)(y)\leq -\frac{3}{4}\zeta(y)+\frac{3\nu^2}{32}I_{B(0,2|\varsigma|+2R)}(y).
\de
Thus, $(\mathbf{H}^{3}_{b_{2}, \s_{2}})$ is right. And it is obvious that $(\mathbf{H}^{4}_{b_{2}, \s_{2}})$ holds.

Finally, following the proof of Theorem \ref{ldptime} (simpler), we obtain that when 
\ce
\lim\limits_{\e\rightarrow0}\frac{\g}{\e}=0,
\de
 the solution $\{X^{\e,\g}_t, \e, \g>0\}$ of the slow part for the system (\ref{exameq}) for any $t\in[0,T]$ satisfies the LDP with the speed $\e^{-1}$ and the good rate function $I$ given by
\ce
I(x; x_0, t)=\inf\left\{\frac{1}{2}\int_0^t|z(s)|^2\dif s: z\in L^2([0,T],\mR^n), s.t. X^z_{x_0}(t)=x\right\},
\de
where $(X^z_{x_0},K^z_{x_0})$ is the unique solution of the following multivalued differential equation
\ce\left\{\begin{array}{l}
\dif X^z_{x_0}(t)\in -\p\psi(X^z_{x_0}(t))\dif t+\bar{\s}_{1}z(t)\dif t,\quad  0\leq t\leq T,\\
X^z_{x_0}(0)=x_0\in\overline{\cD(\p\psi)},
\end{array}
\right.
\de
and $\bar{\s}_{1}=\int_{\mR}\cos(y)\mu(\dif y)$.

We mention that if $\psi=0$, the system (\ref{exameq}) is just a fast mean-reverting stochastic volatility model (cf. \cite{ffk}) and the good rate function $I$ given by
\ce
I(x; x_0, t)=\frac{|x_0-x|^2}{2\bar{\s}^2_{1}t},
\de
which is consistent with the result in \cite[Theorem 2.1]{ffk}.
\ex

\section{Appendix}\label{app}

In this section, we present properties of $\bar{b}_1$ and $\bar{a}_1$ and a comparison principle for a general Hamilton-Jacobi-Bellman equation.

\subsection{Properties of $\bar{b}_1$ and $\bar{a}_1$}

We first prepare the following result whose proof is similar to that for Lemma 4.1 in \cite{q1}.

\bl\label{frozesti}
Assume that $(\mathbf{H}^{1}_{b_{2}, \s_{2}})$ and $(\mathbf{H}^{2}_{b_{2}, \s_{2}})$ hold. Then we have for $x, x_1, x_2\in\overline{\cD(A)}, y, y_1, y_2\in\mR^m$,
\ce
&&\mE|Y_{t}^{x,y}|^2\leq |y|^2e^{-\a t}+C(1+|x|^2), \quad t\geq 0,\\
&&\mE|Y_{t}^{x_{1},y_1}-Y_{t}^{x_{2},y_2}|^{2}\leq|y_1-y_2|^2e^{-\a t}+\frac{C}{\a}|x_1-x_2|^2, \quad t\geq 0.
\de
\el

By the above lemma and the definition of $\mu^{x}$, it holds that
\ce
&&\int_{\mR^{m}}|y|^{2}\mu^{x}(\dif y)=\int_{\mR^{m}}\mE|Y_{t}^{x,y}|^{2}\mu^{x}(\dif y)\leq \int_{\mR^{m}}\(|y|^{2}e^{-\a t}+C(1+|x|^{2})\)\mu^{x}(\dif y)\no\\
&=& e^{-\a t}\int_{\mR^{m}}|y|^{2}\mu^{x}(\dif y)+C(1+|x|^{2}),
\de
and furthermore
\be
\int_{\mR^{m}}|y|^{2}\mu^{x}(\dif y)\leq C(1+|x|^{2}).
\label{inu2}
\ee

\bl\label{barb1bara1lip}
Assume that $(\mathbf{H}^1_{b_{1}, \s_{1}})$, $(\mathbf{H}^2_{b_{1}, \s_{1}})$, $(\mathbf{H}^{1}_{b_{2}, \s_{2}})$ and $(\mathbf{H}^{2}_{b_{2}, \s_{2}})$ hold. Then it holds that for $x_1, x_2\in\overline{\cD(A)}$,
\ce
|\bar{b}_1(x_1)-\bar{b}_1(x_2)|+\|\bar{a}_1(x_1)-\bar{a}_1(x_2)\|\leq C|x_1-x_2|.
\de
\el
\begin{proof}
Since the proofs for Lipschitz continuity of $\bar{b}_1$ and $\bar{a}_1$ are similar, we only show Lipschitz continuity of $\bar{a}_1$.

By the definition of $\bar{a}_1(x)$, it holds that for $x_1, x_2\in\overline{\cD(A)}$,
\ce
&&\|\bar{a}_1(x_1)-\bar{a}_1(x_2)\|\\
&=&\left\|\int_{\mR^m}\sigma_1\sigma^*_1(x_1,y)\mu^{x_1}(\dif y)-\int_{\mR^m}\sigma_1\sigma^*_1(x_2,y)\mu^{x_2}(\dif y)\right\|\\
&=&\left\|\lim_{S\rightarrow\infty}\frac{1}{S}\int_0^S\mE\sigma_1\sigma^*_1(x_1,Y_t^{x_1,y_0})\dif t-\lim_{S\rightarrow\infty}\frac{1}{S}\int_0^S\mE\sigma_1\sigma^*_1(x_2,Y_t^{x_2,y_0})\dif t\right\|\\
&\leq&\lim_{S\rightarrow\infty}\frac{1}{S}\int_0^S\mE\|\sigma_1\sigma^*_1(x_1,Y_t^{x_1,y_0})-\sigma_1\sigma^*_1(x_2,Y_t^{x_2,y_0})\|\dif t\\
&\leq&\lim_{S\rightarrow\infty}\frac{1}{S}\int_0^SC\left(|x_1-x_2|+\mE|Y_t^{x_1,y_0}-Y_t^{x_2,y_0}|\right)\dif t\\
&\leq&C|x_1-x_2|+\lim_{S\rightarrow\infty}\frac{1}{S}\int_0^S\left(\mE|Y_t^{x_1,y_0}-Y_t^{x_2,y_0})|^2\right)^{1/2}\dif t\\
&\leq&(C+\frac{C}{\sqrt \a})|x_1-x_2|,
\de
where we use ($\mathbf{H}^2_{b_1, \sigma_1}$) and Lemma \ref{frozesti}. The proof is complete.
\end{proof}

\subsection{A comparison principle}

Consider the following Hamilton-Jacobi-Bellman equation:
\be\left\{\begin{array}{l} 
\p_t u(t,x)\in F(x,\p_x u(t,x)) -\<A(x),\p_x u(t,x)\>,\\
u(0,x)=h(x), \quad x\in\overline{\cD(A)},
\end{array}
\right.
\label{hjbgene}
\ee
where $F(x,p): \overline{\cD(A)} \times \mathbb{R}^n\to\mR$ is a measurable function. We assume:
\begin{enumerate}[$({\bf A}_1)$]
\item $F(x,p)$ is continuous on $\overline{\cD(A)}\times\mR^n$.
\item For any $R<\infty$ and $|p_1|, |p_2|\leq R$
\ce
|F(x,p_1)-F(x,p_2)|\leq \varpi_R(|p_1-p_2|),
\de
uniformly in $x$, where $\varpi_R$ is a continuous increasing function on $[0,\infty)$ and $\varpi_R(0)=0$.
\end{enumerate}

\bt\label{compprinth}
Suppose that $u, v$ are a viscosity subsolution and a viscosity supersolution for Eq.(\ref{hjbgene}), respectively. If $F(x,p)$ satisfies $({\bf A}_1)$ and $({\bf A}_2)$, $u(t,x)$ or $v(t,x)$ is Lipschitz continuous with respect to $x$, and $u(0,x)\leq v(0,x)$ for any $x\in\overline{\cD(A)}$, then it holds that $u\leq v$ on $[0,T]\times\overline{\cD(A)}$ for any $T>0$.
\et
\begin{proof}
Assume that
\ce
\sup\limits_{(t,x)\in(0,T]\times\overline{\cD(A)}}\{u(t,x)-v(t,x)\}>0.
\de
So, we may choose $M>0$ such that 
\be
\sup\limits_{(t,x)\in(0,T]\times\overline{\cD(A)}}\{u_M(t,x)-v_M(t,x)\}>0,
\label{umvmsup}
\ee
where 
$$
u_M(t,x)=u(t,x)-\frac{M}{T-t}, \quad v_M(t,x)=v(t,x)+\frac{M}{T-t}.
$$
Moreover, $u_M$ is a viscosity subsolution of 
\be
\p_t u(t,x)\in F(x,\p_x u(t,x))-\frac{M}{T^2}-\<A(x),\p_x u(t,x)\>,
\label{hjbnet2}
\ee
and $v_M$ is a viscosity supersolution of 
\be
\p_t u(t,x)\in F(x,\p_x u(t,x))+\frac{M}{T^2}-\<A(x),\p_x u(t,x)\>.
\label{hjbpot2}
\ee

For $\varrho,\eta>1$ and $0<\d<1$, define
\ce
\Phi_{\varrho,\d,\eta}(t,x,s,y)=u_M(t,x)-v_M(s,y)-\frac{\varrho}{2}|x-y|^2-\frac{\d}{2}(|x|^2+|y|^2)-\frac{\eta}{2}(t-s)^2.
\de
Note that the function 
$$
u_M(t,x)-v_M(s,y)-\frac{\varrho}{2}|x-y|^2-\frac{\eta}{2}(t-s)^2
$$
is upper semicontinuous and bounded on $((0,T]\times\overline{\cD(A)})^2$. Thus, the standard optimization technique implies that for any $\varrho,\eta>1$ and $\d\in(0,1)$, there exists a $(\hat t_{\varrho,\d,\eta},\hat x_{\varrho,\d,\eta},\hat s_{\varrho,\d,\eta},\\ \hat y_{\varrho,\d,\eta})\in((0,T]\times\overline{\cD(A)})^2$ such that
\ce
N_{\varrho,\d,\eta}:=\Phi_{\varrho,\d,\eta}(\hat t_{\varrho,\d,\eta},\hat x_{\varrho,\d,\eta},\hat s_{\varrho,\d,\eta},\hat y_{\varrho,\d,\eta})=\sup\limits_{((0,T]\times\overline{\cD(A)})^2}\Phi(t,x,s,y).
\de
Moreover, by (\ref{umvmsup}) we know that
\ce
0<N_{\varrho,\d,\eta}<\infty,
\de
and
\ce
N_{\frac \varrho 2,\frac \d 2,\frac \eta 2}&\geq& \Phi_{\frac \varrho 2,\frac \d 2,\frac \eta 2}(\hat t_{\varrho,\d,\eta},\hat x_{\varrho,\d,\eta},\hat s_{\varrho,\d,\eta},\hat y_{\varrho,\d,\eta})\\
&\geq& N_{\varrho,\d,\eta}+ \frac\varrho 4|\hat x_{\varrho,\d,\eta}-\hat y_{\varrho,\d,\eta}|^2+\frac\d 4(|\hat x_{\varrho,\d,\eta}|^2+|\hat y_{\varrho,\d,\eta}|^2)+\frac\eta 4(\hat t_{\varrho,\d,\eta}-\hat s_{\varrho,\d,\eta})^2.
\de
The above inequality implies that
\be
\lim\limits_{\varrho\rightarrow \infty}\lim\limits_{\d\rightarrow 0}\lim\limits_{\eta\rightarrow\infty}\left\{\varrho|\hat x_{\varrho,\d,\eta}-\hat y_{\varrho,\d,\eta}|^2+\d(|\hat x_{\varrho,\d,\eta}|^2+|\hat y_{\varrho,\d,\eta}|^2)+\eta(\hat t_{\varrho,\d,\eta}-\hat s_{\varrho,\d,\eta})^2\right\}=0.
\label{thrlim}
\ee
In terms of (\ref{thrlim}), it holds that
\be
\varrho|\hat x_{\varrho,\d,\eta}-\hat y_{\varrho,\d,\eta}|^2+\d(|\hat x_{\varrho,\d,\eta}|^2+|\hat y_{\varrho,\d,\eta}|^2)+\eta(\hat t_{\varrho,\d,\eta}-\hat s_{\varrho,\d,\eta})^2\leq \frac{1}{2}.
\label{hatxyts}
\ee
Besides, we mention that 
\ce
\Phi_{\varrho,\d,\eta}(\hat t_{\varrho,\d,\eta},\hat x_{\varrho,\d,\eta},\hat s_{\varrho,\d,\eta},\hat y_{\varrho,\d,\eta})\geq \Phi_{\varrho,\d,\eta}(\hat t_{\varrho,\d,\eta},\hat x_{\varrho,\d,\eta},\hat s_{\varrho,\d,\eta},\hat x_{\varrho,\d,\eta}).
\de
So, simple calculation yields that
\ce
\frac \varrho 2|\hat x_{\varrho,\d,\eta}-\hat y_{\varrho,\d,\eta}|^2&\leq& v(\hat s_{\varrho,\d,\eta},\hat x_{\varrho,\d,\eta})-v(\hat s_{\varrho,\d,\eta},\hat y_{\varrho,\d,\eta})+\frac \d 2(|\hat x_{\varrho,\d,\eta}|^2-|\hat y_{\varrho,\d,\eta}|^2)\\
&\leq&L_v|\hat x_{\varrho,\d,\eta}-\hat y_{\varrho,\d,\eta}|+\frac \d 2(|\hat x_{\varrho,\d,\eta}|+|\hat y_{\varrho,\d,\eta}|)|\hat x_{\varrho,\d,\eta}-\hat y_{\varrho,\d,\eta}|\\
&=&\(L_v+\frac \d 2(|\hat x_{\varrho,\d,\eta}|+|\hat y_{\varrho,\d,\eta}|)\)|\hat x_{\varrho,\d,\eta}-\hat y_{\varrho,\d,\eta}|,
\de
and
\be
\frac \varrho 2|\hat x_{\varrho,\d,\eta}-\hat y_{\varrho,\d,\eta}|\leq \(L_v+\frac \d 2(|\hat x_{\varrho,\d,\eta}|+|\hat y_{\varrho,\d,\eta}|)\),
\label{ahatxhatyboun}
\ee
where $L_v$ is the Lipschitz constant for $v$.

Without loss of generality, we assume that $\hat t_{\varrho,\d,\eta}, \hat s_{\varrho,\d,\eta}\in (0,T)$. Let 
\ce
u_1(t,x):=v_M(\hat s_{\varrho,\d,\eta},\hat y_{\varrho,\d,\eta})+\frac{\varrho}{2}|x-\hat y_{\varrho,\d,\eta}|^2+\frac{\d}{2}(|x|^2+|\hat y_{\varrho,\d,\eta}|^2)+\frac{\eta}{2}(t-\hat s_{\varrho,\d,\eta})^2,
\de
and
\ce
u_2(s,y):=u_M(\hat t_{\varrho,\d,\eta},\hat x_{\varrho,\d,\eta})-\frac{\varrho}{2}|\hat x_{\varrho,\d,\eta}-y|^2-\frac{\d}{2}(|\hat x_{\varrho,\d,\eta}|^2+|y|^2)-\frac{\eta}{2}(\hat t_{\varrho,\d,\eta}-s)^2.
\de
Then $(\hat t_{\varrho,\d,\eta},\hat x_{\varrho,\d,\eta})$ is a point of maximum of $u_M(t,x)-u_1(t,x)$ and $(\hat s_{\varrho,\d,\eta},\hat y_{\varrho,\d,\eta})$ is a point of minimum of $v_M(s,y)-u_2(s,y)$. Since $u_M$ is a viscosity subsolution of Eq.(\ref{hjbnet2})  and $v_M$ is a viscosity supersolution of Eq.(\ref{hjbpot2}), we obtain that
\be
\eta(\hat t_{\varrho,\d,\eta}-\hat s_{\varrho,\d,\eta})&\leq& F(\hat x_{\varrho,\d,\eta}, \varrho(\hat x_{\varrho,\d,\eta}-\hat y_{\varrho,\d,\eta})+\d \hat x_{\varrho,\d,\eta})-\frac{M}{2T^2}\no\\
&&-A_*(\hat x_{\varrho,\d,\eta},  \varrho(\hat x_{\varrho,\d,\eta}-\hat y_{\varrho,\d,\eta})+\d \hat x_{\varrho,\d,\eta}), \label{u1subsol}\\
\eta(\hat t_{\varrho,\d,\eta}-\hat s_{\varrho,\d,\eta})&\geq& F(\hat y_{\varrho,\d,\eta},\varrho(\hat x_{\varrho,\d,\eta}-\hat y_{\varrho,\d,\eta})-\d\hat y_{\varrho,\d,\eta})+\frac{M}{2T^2}\no\\
&&-A^*(\hat y_{\varrho,\d,\eta},\varrho(\hat x_{\varrho,\d,\eta}-\hat y_{\varrho,\d,\eta})-\d\hat y_{\varrho,\d,\eta}).\label{u2supsol}
\ee
By subtracting (\ref{u2supsol}) from (\ref{u1subsol}), it holds that
\be
\frac{M}{T^2}&\leq& F(\hat x_{\varrho,\d,\eta}, \varrho(\hat x_{\varrho,\d,\eta}-\hat y_{\varrho,\d,\eta})+\d \hat x_{\varrho,\d,\eta})-F(\hat y_{\varrho,\d,\eta},\varrho(\hat x_{\varrho,\d,\eta}-\hat y_{\varrho,\d,\eta})-\d\hat y_{\varrho,\d,\eta})\no\\
&&-\(A_*(\hat x_{\varrho,\d,\eta},  \varrho(\hat x_{\varrho,\d,\eta}-\hat y_{\varrho,\d,\eta})+\d \hat x_{\varrho,\d,\eta})-A^*(\hat y_{\varrho,\d,\eta},\varrho(\hat x_{\varrho,\d,\eta}-\hat y_{\varrho,\d,\eta})-\d\hat y_{\varrho,\d,\eta})\).\no\\
\label{FAes}
\ee

In the following, we observe
$$
A_*(\hat x_{\varrho,\d,\eta},  \varrho(\hat x_{\varrho,\d,\eta}-\hat y_{\varrho,\d,\eta})+\d \hat x_{\varrho,\d,\eta})-A^*(\hat y_{\varrho,\d,\eta},\varrho(\hat x_{\varrho,\d,\eta}-\hat y_{\varrho,\d,\eta})-\d\hat y_{\varrho,\d,\eta}).
$$
On the one hand, by Lemma 2 in \cite{Za}, it holds that
\ce
&&A_*(\hat x_{\varrho,\d,\eta},\varrho(\hat x_{\varrho,\d,\eta}-\hat y_{\varrho,\d,\eta})+\d \hat x_{\varrho,\d,\eta})=\inf\limits_{x^*\in A(\hat x_{\varrho,\d,\eta})}\<x^*,\varrho(\hat x_{\varrho,\d,\eta}-\hat y_{\varrho,\d,\eta})+\d \hat x_{\varrho,\d,\eta}\>,\\
&&A^*(\hat y_{\varrho,\d,\eta},\varrho(\hat x_{\varrho,\d,\eta}-\hat y_{\varrho,\d,\eta})-\d\hat y_{\varrho,\d,\eta})=\sup\limits_{y^*\in A(\hat y_{\varrho,\d,\eta})}\<y^*,\varrho(\hat x_{\varrho,\d,\eta}-\hat y_{\varrho,\d,\eta})-\d\hat y_{\varrho,\d,\eta}\>.
\de
On the other hand, without restricting the generality, by $({\bf H}_A)$ we can suppose that $0\in A(0)$. So, for any $x^*\in A(\hat x_{\varrho,\d,\eta})$ and $y^*\in A(\hat y_{\varrho,\d,\eta})$, by Lemma \ref{equi},
\ce
\<x^*,\varrho(\hat x_{\varrho,\d,\eta}-\hat y_{\varrho,\d,\eta})+\d \hat x_{\varrho,\d,\eta}\>&\geq& \<x^*,\varrho(\hat x_{\varrho,\d,\eta}-\hat y_{\varrho,\d,\eta})\>\geq \<y^*,\varrho(\hat x_{\varrho,\d,\eta}-\hat y_{\varrho,\d,\eta})\>\\
&\geq&\<y^*,\varrho(\hat x_{\varrho,\d,\eta}-\hat y_{\varrho,\d,\eta})-\d\hat y_{\varrho,\d,\eta}\>.
\de
The above deduction implies that
\be
A_*(\hat x_{\varrho,\d,\eta},  \varrho(\hat x_{\varrho,\d,\eta}-\hat y_{\varrho,\d,\eta})+\d \hat x_{\varrho,\d,\eta})-A^*(\hat y_{\varrho,\d,\eta},\varrho(\hat x_{\varrho,\d,\eta}-\hat y_{\varrho,\d,\eta})-\d\hat y_{\varrho,\d,\eta})\geq 0.
\label{a*a*sign}
\ee

Finally, combining (\ref{FAes}) with (\ref{a*a*sign}), we conclude that
\be
\frac{M}{T^2}\leq F(\hat x_{\varrho,\d,\eta}, \varrho(\hat x_{\varrho,\d,\eta}-\hat y_{\varrho,\d,\eta})+\d \hat x_{\varrho,\d,\eta})-F(\hat y_{\varrho,\d,\eta},\varrho(\hat x_{\varrho,\d,\eta}-\hat y_{\varrho,\d,\eta})-\d\hat y_{\varrho,\d,\eta}).
\label{mfes}
\ee
(\ref{hatxyts}) and (\ref{ahatxhatyboun}) infer that
\ce
&&\lim\limits_{\varrho\rightarrow\infty}\hat x_{\varrho,\d,\eta}=\lim\limits_{\varrho\rightarrow\infty}\hat y_{\varrho,\d,\eta}=:l_{\d,\eta},\\
&&\lim\limits_{\varrho\rightarrow\infty}\varrho(\hat x_{\varrho,\d,\eta}-\hat y_{\varrho,\d,\eta})=:L_{\d,\eta}.
\de
Taking the limit on two sides of (\ref{mfes}) as $\varrho\rightarrow\infty$, we obtain that
\ce
\frac{M}{T^2}\leq F(l_{\d,\eta},L_{\d,\eta}+\d l_{\d,\eta})-F(l_{\d,\eta},L_{\d,\eta}-\d l_{\d,\eta}).
\de
Again from (\ref{hatxyts}) and (\ref{ahatxhatyboun}), it follows that 
$$
\sup\limits_{\d}|L_{\d,\eta}+\d l_{\d,\eta}|\leq C, \quad \sup\limits_{\d}|L_{\d,\eta}-\d l_{\d,\eta}|\leq C.
$$
Thus, by $({\bf A}_2)$ and (\ref{hatxyts}), we have that
\ce
\frac{M}{T^2}\leq \varpi_C(2\d |l_{\d,\eta}|)\leq \varpi_C(\sqrt{2\d}).
\de
Letting $\d\rightarrow 0$, by the properties of $\varpi_C$, one can get that
$$
\frac{M}{T^2}\leq0,
$$
which contradicts with the fact that $M>0$. Therefore, we draw the conclusion that $u\leq v$ on $[0,T]\times\overline{\cD(A)}$ for any $T>0$. The proof is complete.
\end{proof}

\end{document}